\documentclass{amsproc}

\usepackage[swedish,english]{babel}
\usepackage[utf8]{inputenc}
\usepackage[T1]{fontenc}
\usepackage{amsmath,amssymb,amsfonts}
\usepackage{url}
\usepackage[dvipsnames,usenames]{color}
\usepackage{array}
\usepackage[pdftex]{graphicx}
\usepackage{tikz}
\usetikzlibrary{calc,decorations.markings} 
\usepackage{accents}
\usepackage{tabu}
\usepackage[normalem]{ulem} 
\usepackage{mathdots}

\usepackage{tikz}
\usetikzlibrary{arrows,positioning,matrix} 

\newcolumntype{L}{>{\displaystyle}l}
\newcolumntype{C}{>{\displaystyle}c}
\newcolumntype{R}{>{\displaystyle}r}

\renewcommand{\tfrac}{\genfrac{}{}{}1}

\newcommand{\R}{\mathbb R}
\newcommand{\N}{\mathbb N}
\newcommand{\Z}{\mathbb Z}
\newcommand{\C}{\mathbb C}

\newcommand{\cont}{\mathcal C}
\newcommand{\D}{\mathbb D}
\renewcommand{\Im}{\mathrm{Im}}
\renewcommand{\Re}{\mathrm{Re}}
\def\E{{\mathrm{e}}}

\def\til{\widetilde}
\newcommand{\hol}{\mathcal O}
\newcommand{\pois}{\mathcal P}
\newcommand{\cpois}{\mathcal{Q}}
\newcommand{\Res}{\mathrm{Res}}
\def\I{\mathfrak{i}}
\newcommand{\supp}{\mathrm{supp}}
\newcommand{\diff}{\mathrm{d}}
\renewcommand{\bar}{\overline}
\newcommand{\sh}{\mathrm{sh}}
\newcommand{\ch}{\mathrm{ch}}
\renewcommand{\th}{\mathrm{th}}

\newcommand{\ntto}{\:\scriptsize{\xrightarrow{\vee\:}}\:}
\newcommand{\nntto}{\:\scriptsize{\xrightarrow{\wedge\:}}\:}

\newcommand{\hardy}{\mathrm{H}}

\newcommand{\qRe}[1]{{{#1}_\mathrm{re}}}
\newcommand{\qIm}[1]{{{#1}_\mathrm{im}}}

\newcommand{\ie}{\textit{i.e.}\/ } 
\newcommand{\eg}{\textit{e.g.}\/ } 
\newcommand{\cf}{\textit{cf.}\/ } 

\allowdisplaybreaks

\theoremstyle{definition} 
\newtheorem{define}{Definition}[section]
\newtheorem{example}[define]{Example}
\newtheorem{remark}[define]{Remark}

\theoremstyle{plain} 
\newtheorem{lemma}[define]{Lemma}
\newtheorem{thm}[define]{Theorem}
\newtheorem{prop}[define]{Proposition}
\newtheorem{coro}[define]{Corollary}

\numberwithin{equation}{section}

\begin{document}

\title[On quasi-Herglotz functions in one variable]{On quasi-Herglotz functions in one variable}

\author{Annemarie Luger}
\address{Annemarie Luger, Department of Mathematics, Stockholm University, SE-106 91 Stockholm, Sweden}
\curraddr{}
\email{luger@math.su.se}
\thanks{The authors were partially supported by the Swedish Foundation for Strategic Research, grant nr. AM13-0011.}

\author{Mitja Nedic}
\address{Mitja Nedic, Department of Mathematics, Stockholm University, SE-106 91 Sto\-ckholm, Sweden, orc-id: 0000-0001-7867-5874}
\curraddr{}
\email{mitja@math.su.se}
\thanks{}


\date{2018-12-07} 

\begin{abstract}
In this paper,  the class of (complex) quasi-Herglotz functions is introduced as the complex vector space generated by the convex cone of ordinary Herglotz functions. We prove  characterization theorems, in particular, an analytic characterization. The subclasses of quasi-Herglotz functions that are identically zero in one half-plane as well as rational quasi-Herglotz functions are investigated in detail. Moreover, we relate to other areas such as weighted Hardy spaces, definitizable functions, the Cauchy transform on the unit circle and sum-rule identities.
\end{abstract}

\maketitle

\section{Introduction}\label{sec:introduction}

Holomorphic functions in a domain in $\C$ form the centrepiece of one-dimensional complex analysis. However, in certain applications, we often restrict ourselves to a smaller subclass of functions where we can derive additional information. One such prominent class of functions is the class of \emph{Herglotz functions} (also called \emph{Nevanlinna functions}, \emph{Herglotz-Nevanlinna functions}, \emph{Pick functions}, \emph{R-functions}, etc.), which consists of all holomorphic functions on $\C\setminus\R$ that map the upper and lower half-planes to the closed upper and lower half-planes, respectively, \cf Definition \ref{def:Herglotz}. These functions appear at many places in pure mathematics as well as in several applications. Generally speaking, they are deeply connected with extension theory of symmetric operators \cite{AkhiezerGlazman1993,LangerTextorius1977}. In particular, they appear, for example, as Titchmarsh-Weyl $m$-functions for differential operators, see \eg \cite{EckhardtKostenkoTeschl2017,Everitt1972}, which are also related to scattering problems \cite{AlbeverioHrynivMykytyuk2012} or in connection with the moment problem \cite{Akhiezer1965,DamanikKillipSimon2010}. They (or rational transformations) appear also as certain transfer functions of input/state/output linear systems, see \eg \cite{BallBolotnikov2010} for a discussion that includes also multivariable situations. Moreover, Herglotz functions are also used \eg in number theory \cite{KurlbergUeberschaer2017}. In applications, they play an important role \eg as the Fourier transform of the impulse response of passive systems, which can describe the (one-port) system completely \cite{IvanenkoETAL2019a,Zemanian1963} and are essential \eg in for physical limitations using so-called sum rules \cite{BernlandEtal2011}. Another area of application is within homogenization \cite{Milton2002} to name but a few.

The effectiveness of this particular class of functions can partly be attributed to the classical integral representations theorem, see \eg \cite{Cauer1932,KacKrein1974}, which states that a function $h$ is Herglotz function if and only if it can be written in the form
$$h(z) = a + b\:z + \frac{1}{\pi}\int_\R\frac{1+t\:z}{t-z}\diff\nu(t),$$
for the details see Theorem \ref{thm:intRep_1var} and representation \eqref{eq:intRep_1var_finite}.

In this article, we consider a larger class of functions, namely \emph{quasi-Herglotz functions}. These are introduced as the complex vector space generated by the convex cone of ordinary Herglotz functions, \cf Definition \ref{def:quasi_Herglotz}. They admit an integral representation of the form above where the parameters $a$, $b$ and $\nu$ not only lack a sign constraint as in Theorem \ref{thm:intRep_1var}, but are allowed to be complex, \cf Theorem \ref{thm:intRep_1var_quasi}. We also give different characterizations, discuss some important subclasses and, finally, relate our results to other areas.

Although (real or complex) linear combinations of ordinary Herglotz functions appear at several places, they have not been investigated systematically so far. They appear, for example, in connection with non-passive and active systems, \eg active exterior cloaking or non-passive gain media \cite{GuevaraVasquezMiltonOnofrei2012,IvanenkoETAL2019b}. Particularly in \cite{IvanenkoETAL2019b}, it is of special interest that quasi-Herglotz functions, like ordinary Herglotz functions, allow for the application of convex optimization methods in electrostatic modeling. Another example is Kre\u{\i}n's spectral shift function \cite{BehrndtGesztesyNakamura2018,GesztesyMakarovNaboko1999,Krein1953}, which can also be viewed in the context of quasi-Herglotz functions.

The structure of the paper is as follows. After the introduction in Section \ref{sec:introduction}, we give the necessary background information about ordinary Herglotz functions in Section \ref{sec:background}. In Section \ref{sec:quasi_HG}, we introduce quasi-Herglotz functions, collect some basic properties and prove an integral representation. Then, in Section \ref{sec:analytic}, we answer the question which holomorphic functions are actually quasi-Herglotz functions by an analytic characterizarion, \cf Theorem \ref{thm:quasi_Tumarkin_Vladimirov}. In Section \ref{subsec:sym_and_uni},  symmetry properties are discussed as well as to what extent are quasi-Herglotz functions determined by their values in only one half-plane. Section \ref{sec:identically_zero} focuses on quasi-Herglotz functions that are identically zero in one half-plane. Rational quasi-Herglotz functions are discussed in Section \ref{sec:rational} and Section \ref{sec:comparisons} relates quasi-Herglotz functions to other areas of analysis. In particular, Section \ref{subsec:Hardy} highlights the intersection between quasi-Herglotz functions that are identically zero in one half-plane and a weighted Hardy space $\hardy^1$, in Section \ref{subsec:definitizable} the connection with definitizable functions is explained, and Section \ref{subsec:Cauchy} presents how the Cauchy transform of a complex Borel measure on the unit circle $S^1$ may be viewed as a special case of a quasi-Herglotz function, \cf Theorem \ref{thm:Cauchy}. Finally, in Section \ref{subsec:sumrules}, we discuss the existence of sum-rules for real quasi-Herglotz functions.

\section{Background}\label{sec:background}

We start be recalling the definition of a Herglotz function, \cf \cite{KacKrein1974}.

\begin{define}\label{def:Herglotz}
A function $h\colon \C\setminus\R \to \C$ is called \emph{Herglotz function} if it is holomorphic with  
$$\frac{\Im[h(z)]}{\Im[z]} \geq 0 \quad\text{and}\quad \bar{h(z)} = h(\bar{z})$$
for all $z \in \C\setminus\R$.
\end{define}

We stress that this definition, where {the function} $h$ is defined both in the upper and the lower half-planes, does not constitute the only possible way  and another - equivalent - way of defining Herglotz functions will be discussed in Section \ref{subsec:sym_and_uni}.

The classical integral representation formula for Herglotz functions may be presented in the following formulation, \cf \cite{Cauer1932,KacKrein1974}.

\begin{thm}\label{thm:intRep_1var}
A function $h\colon \C\setminus\R \to \C$ is a Herglotz function if and only if {$h$ can be written}, for every $z \in \C\setminus\R$, {as}
\begin{equation}
\label{eq:intRep_1var}
h(z) = a + b\:z + \frac{1}{\pi}\int_\R K(z,t)\diff\mu(t),
\end{equation}
where the kernel $K\colon (\C\setminus\R) \times \R \to \C$ is defined as
$$K(z,t) := \frac{1}{t-z} - \frac{t}{1+t^2} = \frac{1+t\:z}{(t-z)(1+t^2)}$$
and $a \in \R$, $b \geq 0$ and $\mu$ is a positive Borel measure on $\R$ satisfying the growth condition
\begin{equation}
\label{eq:growth_1var}
\int_\R\frac{1}{1+t^2}\diff\mu(t) < \infty.
\end{equation}
\end{thm}

Furthermore, representation \eqref{eq:intRep_1var} is unique for a given function $h$ and the collection $(a,b,\mu)$ of the representing parameters is called the \emph{data} corresponding to the function $h$ in the sense of representation \eqref{eq:intRep_1var}. Moreover, given a function $h$, its data can be obtained in the following way. It holds that
$$a = \Re[h(\I)]$$
and that the measure $\mu$ is given by the Stieltjes inversion formula \cite{KacKrein1974}, \ie
\begin{equation}
\label{eq:Stieltjes_inversion_herglotz}
\lim\limits_{y \to 0^+}\int_\R g(x)\Im[h(x+\I\:y)]\diff x = \int_\R g(t)\diff\mu(t)
\end{equation}
for any $\cont^1$-function $g\colon \R \to \R$ such that there exists a constant $C \in \R$ with the property that $|g(x)| \leq C(1+x^2)^{-1}$ for all $x \in \R$.

\begin{figure}[!ht]
\centering
\begin{tikzpicture}[scale=1.1]
\fill[fill=black!10!white] (-1.5,2.179) -- (1.5,0) -- (4.5,2.179);
\fill[fill=black!10!white] (-2.2,-2.179) -- (-0.5,0) -- (1.2,-2.179);
\draw[help lines,->] (-3.2,0) -- (4.7,0) node[above] {$x$};
\draw[help lines,->] (0,-2.5) -- (0,2.5) node[right] {$\I\,y$};
\draw[-] (-1.5,2.179) -- (1.5,0) -- (4.5,2.179);
\draw[-] (-2.2,-2.179) -- (-0.5,0) -- (1.2,-2.179);
\draw (1.5,0) node[below] {$t_1$};
\draw (-0.5,0) node[above] {$t_2$};
\draw [dashed,domain=0:36] plot ({1.5+0.75*cos(\x)}, {0.75*sin(\x)});
\draw [dashed,domain=0:51] plot ({-0.5-0.75*cos(\x)}, {-0.75*sin(\x)});
\draw (2.5,0.3) node {$\theta_1$};
\draw (-1.5,-0.35) node {$\theta_2$};
\end{tikzpicture}
\caption{An upper Stoltz domain with centre $t_1$ and angle $\theta_1$ and a lower Stoltz domain with centre $t_2$ and angle $\theta_2$.}
\label{fig:Stoltz_domain}
\end{figure}
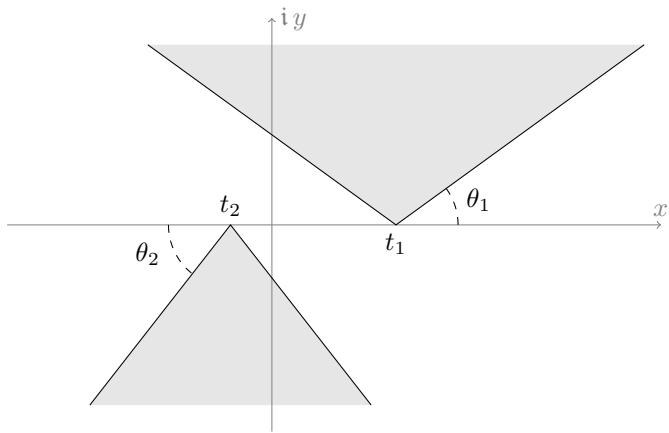

To be able to describe the parameter $b$, we recall first the definition of  a non-tangential limit. An \emph{upper Stoltz domain} with centre $t_0 \in \R$ and angle $\theta \in (0,\frac{\pi}{2}]$ is the set
$$\{z \in \C^+~|~\theta \leq \arg(z-t_0) \leq \pi-\theta\},$$
see Figure \ref{fig:Stoltz_domain}. The symbol $z \ntto t_0$ then denotes the limit $z \to t_0$ in any upper Stoltz domain with centre $t_0$ and the symbol $z \ntto \infty$ denotes the limit $|z| \to \infty$ in any upper Stoltz domain with centre $0$. \emph{Lower Stoltz domains} are defined analogously and non-tangential limits in a lower Stoltz domain are denoted by $z \nntto t_0$ or $z \nntto \infty$. Note that in the literature, slightly different notations are also used for these non-tangential limits. Examples of upper and lower Stoltz domains are visualized in Figure \ref{fig:Stoltz_domain}.

We may now obtain the constant $b$ from representation \eqref{eq:intRep_1var} as
\begin{equation}\label{eq:b_parameter_Herglotz}
b = \lim\limits_{z \ntto \infty}\frac{h(z)}{z} = \lim\limits_{z \nntto \infty}\frac{h(z)}{z}.
\end{equation}

Sometimes, it is more convenient to rewrite representation \eqref{eq:intRep_1var} in such a way that the measure $\mu$ is replaced by a finite measure. This can be done, for example, by defining, for any measure $\mu$ as before, a measure $\nu$ via
\begin{equation}\label{eq:measure_make_finite}
\diff\nu(t) := \frac{1}{1+t^2}\diff\mu(t).
\end{equation}
Representation \eqref{eq:intRep_1var} may thus be rewritten as
\begin{equation}
\label{eq:intRep_1var_finite}
h(z) = a + b\:z + \frac{1}{\pi}\int_\R \til{K}(z,t)\diff\nu(t),
\end{equation}
where the kernel $\til{K}\colon (\C\setminus\R) \times \R \to \C$ is defined by
$$\til{K}(z,t) := \frac{1+t\:z}{t-z} = (1+t^2)K(z,t).$$

\section{Quasi-Herglotz functions and basic properties}\label{sec:quasi_HG}

The set of ordinary Herglotz function is a convex cone, \ie any non-negative linear combination of Herglotz functions is again a Herglotz function. The introduction of quasi-Herglotz function extends the set of Herglotz functions to a complex vector space.

\begin{define}\label{def:quasi_Herglotz}
A function $q\colon \C\setminus\R \to \C$ is called a \emph{quasi-Herglotz function} if there exist Herglotz functions $h_1,h_2,h_3$ and $h_4$, such that it holds, for all $z \in \C\setminus\R$, that
\begin{equation}
    \label{eq:quasi_combination}
    q(z) = (h_1-h_2)(z) + \I\,(h_3-h_4)(z).
\end{equation}
\end{define}

\begin{remark}
Decomposition \eqref{eq:quasi_combination} of a function $q$ is not unique, as one may add a fixed Herglotz function $h$ to both functions $h_1$ and $h_2$ (or to $h_3$ and $h_4$) and still get the same function $q$ as a result.
\end{remark}

When considering only real-linear combinations of ordinary Herglotz functions, the following subclass arises, \cf \cite[Def. 2.1]{IvanenkoETAL2019b}.

\begin{define}\label{def:real_quasi_Herglotz}
A quasi-Herglotz function is called \emph{real} if for any four function $h_1,h_2,h_3$ and $h_4$ satisfying formula \eqref{eq:quasi_combination} it holds that $h_3 = h_4$.
\end{define}

\subsection{An integral representation for quasi-Herglotz functions}

Our first result is an integral representation theorem for quasi-Herglotz functions which, in parts, is a straight forward consequence of Theorem \ref{thm:intRep_1var}. The crucial point in the statement, however, concerns the uniqueness of the data. Moreover, it shows that quasi-Herglotz functions are the largest class of holomorphic functions on $\C\setminus\R$ admitting an integral representation of the same form as ordinary Herglotz functions.

\begin{thm}\label{thm:intRep_1var_quasi}
A function $q\colon\C\setminus\R \to \C$ is a quasi-Herglotz function if and only if $q$ can be written, for every $z \in \C\setminus\R$, as
\begin{equation}\label{eq:intRep_1var_quasi}
q(z) = a+b\:z + \frac{1}{\pi}\int_\R \til{K}(z,t)\diff\nu(t),
\end{equation}
where $a,b \in \C$ and $\nu$ is a complex Borel measure on $\R$. Furthermore, this representation is unique, \ie each quasi-Herglotz function $q$ is uniquely determined by its data-triple $(a,b,\nu)$.
\end{thm}

\begin{remark}
A complex (or signed) Borel measure is finite by definition.
\end{remark}

\proof
First, if $q$ is a quasi-Herglotz function, then it can be written as
$$q(z) = (h_1-h_2)(z) + \I\,(h_3-h_4)(z),$$
where $h_1,h_2,h_3,h_4$ are four Herglotz functions. Each of these four Herglotz functions admits an integral representation of the form \eqref{eq:intRep_1var_finite}, and combining these gives an integral representation of the form \eqref{eq:intRep_1var_quasi} for the function $q$.

Conversely, any function $q$, defined by the integral representation \eqref{eq:intRep_1var_quasi} on $\C\setminus\R$ may be written as a combination of four Herglotz functions by splitting $a = a_1-a_2 + \I\,(a_3-a_4) \in \C$, $b = b_1-b_2 + \I\,(b_3-b_4) \in \C$ and $\nu = \nu_1 - \nu_2 + \I\,(\nu_3-\nu_4)$, where $a_j,b_j \geq 0$ and $\nu_j$ are finite positive Borel measure. The Herglotz function $h_j$ can then be taken as given by the data $(a_j,b_j,\nu_j)$ in the sense of representation \eqref{eq:intRep_1var_finite}.

Therefore, it remains to show that the data corresponding to a function $q$ is uniquely determined by the function. To that end, suppose that there exists two sets of data, namely $(a,b,\nu)$ and $(a',b',\nu')$, such that the function $q$ admits a representation of the form \eqref{eq:intRep_1var_quasi} with respect to both sets of data.

If this is the case, then considering the limit
$$\lim\limits_{z \ntto \infty}\frac{q(z)}{z}$$
via representation \eqref{eq:intRep_1var_quasi} using the two data sets $(a,b,\nu)$ and $(a',b',\nu')$ yields that $b = b'$. Considering the expression $q(\I) + q(-\I)$ in an analogous way yields $a = a'$.

Thus, it remains to show that $\nu = \nu'$. To do this, it suffices to show that a complex Borel measure $\eta$ on $\R$ is identically zero whenever
\begin{equation}\label{eq:eta_1}
\int_\R \til{K}(z,t)\diff\eta(t) = 0
\end{equation}
for all $z \in \C\setminus\R$. Investigating the conjugate of equality \eqref{eq:eta_1} yields
$$0 = \int_\R \bar{\til{K}(z,t)}\diff\bar{\eta}(t) = \int_\R \til{K}(\bar{z},t)\diff\bar{\eta}(t) = \int_\R \til{K}(\zeta,t)\diff\bar{\eta}(t)$$
where $\zeta \in \C\setminus\R$, $\bar{\eta} := \qRe{\eta} - \I\,\qIm{\eta}$ and $\qRe{\eta}$ and $\qIm{\eta}$ are two signed measures on $\R$ such that $\eta = \qRe{\eta} + \I\,\qIm{\eta}$. Therefore, if a complex measure $\eta$ satisfies equality \eqref{eq:eta_1}, so does its conjugate measure $\bar{\eta}$. Hence, it follows that
$$\int_\R \til{K}(z,t)\diff\qRe{\eta}(t) = \int_\R \til{K}(z,t)\diff\qIm{\eta}(t) = 0$$
for all $z \in \C\setminus\R$.

As such, we may assume, without loss of generality, that $\eta$ is a signed measure. Considering, again, the conjugate of equality \eqref{eq:eta_1}, we infer that
\begin{multline}\label{eq:eta_2}
    0 = \int_\R \til{K}(z,t)\diff\eta(t) - \bar{\int_\R \til{K}(z,t)\diff\eta(t)} \\
    = \int_\R \til{K}(z,t)\diff\eta(t) - \int_\R \til{K}(\bar{z},t)\diff\eta(t) = 2\I\int_\R(1+t^2)\pois(z,t)\diff\eta(t),
\end{multline}
where $\pois$ denotes the Poisson kernel of the upper half-plane $\C^+$, \ie
\begin{equation}
\label{eq:poisson_kernels}
\pois(z,t) := \frac{\Im[z]}{|t-z|^2}.
\end{equation}

Let $\eta_1$ and $\eta_2$ be two finite positive Borel measures on $\R$, such that $\eta = \eta_1-\eta_2$, and define, further, two positive Borel measures $\rho_1$ and $\rho_2$ on $\R$ by setting
$$\diff\rho_j(t) := (1+t^2)\diff\eta_j(t)$$
for $j=1,2$. Then, we infer from equation \eqref{eq:eta_2} that
$$\int_\R\pois(z,t)\diff\rho_1(t) = \int_\R\pois(z,t)\diff\rho_2(t)$$
for all $z \in \C\setminus\R$. An elementary property of the Poisson kernel, see \eg \cite[p. 111]{Koosis1998}, implies now that $\rho_1 \equiv \rho_2$, yielding back that $\eta_1 \equiv \eta_2$ or, in other words, that $\eta \equiv 0$. This finishes the proof.
\endproof

\begin{remark}
Note that in order to show the uniqueness of the representation, values of $q$ in both the upper half plane $\C^+$ and the lower half plane $\C^-$ have been used. 
\end{remark}

The following corollary is now an immediate consequence of the preceding proof.

\begin{coro}\label{coro:a_and_b_parameters}
The numbers $a$ and $b$ from Theorem \ref{thm:intRep_1var_quasi} are equal to
\begin{equation}\label{eq:a_constant}
    a = \frac{1}{2}\big(q(\I) + q(-\I)\big)
\end{equation}
and
\begin{equation}\label{eq:b_constant}
    b = \lim\limits_{z \ntto \infty}\frac{q(z)}{z} = \lim\limits_{z \nntto \infty}\frac{q(z)}{z}.
\end{equation}
\end{coro}

It turns out that the measure $\nu$ also satisfies an inversion formula similar to formula \eqref{eq:Stieltjes_inversion_herglotz}. In fact, we are going to prove later in Proposition \ref{prop:Stieltjes_inversion} that
$$\lim\limits_{y \to 0^+}\int_\R g(x)\tfrac{1}{2\I}(q(x+\I\:y)-q(x-\I\:y))\diff x = \int_\R g(t)(1+t^2)\diff\nu(t)$$
for all admissible functions $g$. However, the proof of this result requires {the introduction of some additional concepts in Section \ref{subsec:Stieltjes}.}

For real quasi-Herglotz functions, Theorem \ref{thm:intRep_1var_quasi} provides us with the following additional corollaries.

\begin{coro}\label{coro:real_characterization}
A quasi-Herglotz function $q$ is real if and only if its representing parameters are real, \ie $a,b \in \R$ and $\nu$ is a signed Borel measure on $\R$.
\end{coro}

\begin{coro}
Let $q$ be a real quasi-Herglotz functions. Then it holds, for every $z \in \C^+$, that
\begin{equation*}
\label{eq:intRep_re_part_real_quasi}
\Re[q(z)] = a + b\:\Re[z] + \frac{1}{\pi}\int_\R\left((1+t^2)\cpois(z,t) - t\right)\diff\nu(t)
\end{equation*}
and
\begin{equation*}
\label{eq:intRep_im_part_real_quasi}
\Im[q(z)] = b\:\Im[z] + \frac{1}{\pi}\int_\R(1+t^2)\pois(z,t)\diff\nu(t),
\end{equation*}
where $\pois$ and $\cpois$ denote the Poisson kernel and conjugate Poisson kernel of the upper half-plane in one variable, respectively, \ie
$$\pois(z,t) := \frac{\Im[z]}{|t-z|^2} \quad\text{and}\quad \cpois(z,t) := \frac{\Re[z]+t}{|t-z|^2},$$
where $z \in \C\setminus\R$ and $t \in \R$.
\end{coro}

\subsection{Stieltjes inversion formula}\label{subsec:Stieltjes}

For ordinary Herglotz functions, it is the imaginary part of the function that determines the representing measures as evident from formula \eqref{eq:Stieltjes_inversion_herglotz}. For quasi-Herglotz functions, it's imaginary part does no longer play the same role. Instead, we need an appropriate substitute and to that end, we consider the following definitions.

\begin{define}\label{def:conjugate_function}
Let $q$ be a quasi-Herglotz function. Then, its \emph{conjugate function} $\bar{q}$ is the quasi-Herglotz function given by the conjugate parameters of $q$, \ie if $q$ is represented by the data $(a,b,\nu)$ in the sense of Theorem \ref{thm:intRep_1var_quasi}, then $\bar{q}$ represented by the data $(\bar{a},\bar{b},\bar{\nu})$.
\end{define}

\begin{define}
Let $q$ be a quasi-Herglotz function. Then, its \emph{quasi-real part} $\qRe{q}$ and \emph{quasi-imaginary part} $\qIm{q}$ are defined, {for $z \in \C\setminus\R$}, as
$$\qRe{q}(z) := \frac{1}{2}(q(z) + \bar{q}(z)) \quad\text{and}\quad \qIm{q}(z) := \frac{1}{2\I}(q(z)-\bar{q}(z)).$$
\end{define}

We observe that if the function $q$ is represented by the data $(a,b,\nu)$, then its quasi-real part $\qRe{q}$ is represented by the data $(\Re[a],\Re[b],\qRe{\nu})$ and its quasi-imaginary part $\qIm{q}$ is represented by the data $(\Im[b],\Im[b],\qIm{\nu})$. In other words, the quasi-real and quasi-imaginary parts of a quasi-Herglotz function are real quasi-Herglotz functions. Note, furthermore, that
$$q(z) = \qRe{q}(z) + \I\:\qIm{q}(z)$$
for every $z \in \C\setminus\R$, while
$$\qRe{q}(z) \neq \Re[q(z)] \quad\text{and}\quad \qIm{q}(z) \neq \Im[q(z)]$$
in general.

We will now show that real quasi-Herglotz functions satisfy a direct analogue of the Stieltjes inversion formula \eqref{eq:Stieltjes_inversion_herglotz} for ordinary Herglotz functions, which will then be used to determine the analogue of the Stieltjes inversion formula for (non-real) quasi-Herglotz functions.

\begin{lemma}\label{lem:Stieltjes_inversion_real}
Let $q$ be a real quasi-Herglotz function. Then, its representing measure $\nu$ satisfies 
\begin{equation}
\label{eq:Stieltjes_inversion_quasi_real}
    \lim\limits_{y \to 0^+}\int_\R g(x)\Im[q(x+\I\:y)]\diff x = \int_\R(1+t^2)g(t)\diff\nu(t)
\end{equation}
for any $\cont^1$-function $g\colon \R \to \R$ such that there exists a constant $C \in \R$ with the property that $|g(x)| \leq C(1+x^2)^{-1}$ for all $x \in \R$.
\end{lemma}

The proof of this lemma follows immediately by combining the corresponding results for any two ordinary Herglotz functions $h_1$ and $h_2$ such that $q = h_1 - h_2$.

\begin{prop}\label{prop:Stieltjes_inversion}
The measure $\nu$ from Theorem \ref{thm:intRep_1var_quasi} satisfies the formula
\begin{equation}\label{eq:quasi_stieltjes_v1}
\lim\limits_{y \to 0^+}\int_\R g(x)\tfrac{1}{2\I}(q(x+\I\:y)-q(x-\I\:y))\diff x = \int_\R g(t)(1+t^2)\diff\nu(t),
\end{equation}
where $g\colon \R \to \C$ is a $\cont^1$-function such that there exists a constant $C \geq 0$ so that $|g(x)| \leq C(1+x^2)^{-1}$ for all $x \in \R$.
\end{prop}

\proof
By Lemma \ref{lem:Stieltjes_inversion_real}, we know that the functions $\qRe{q}$ and $\qIm{q}$ satisfy the formulas
\begin{equation}\label{eq:stieltjes_calc1}
\lim\limits_{y \to 0^+}\int_\R g_1(x)\Im[\qRe{q}(x+\I\:y)]\diff x = \int_\R g_1(t)(1+t^2)\diff(\qRe{\nu})(t)
\end{equation}
and
\begin{equation}\label{eq:stieltjes_calc2}
\lim\limits_{y \to 0^+}\int_\R g_2(x)\Im[\qIm{q}(x+\I\:y)]\diff x = \int_\R g_2(t)(1+t^2)\diff(\qIm{\nu})(t),
\end{equation}
where $g_1,g_2\colon \R \to \R$ are two $\cont^1$-functions satisfying the assumption of the Stieltjes inversion formula.

Adding now to formula \eqref{eq:stieltjes_calc1} an $\I$-multiple of formula \eqref{eq:stieltjes_calc2} yields
\begin{multline*}
\lim\limits_{y \to 0^+}\int_\R(g_1 + \I\:g_2)(x)\left(\Im[\qRe{q}(x+\I\:y)]+\I\:\Im[\qIm{q}(x+\I\:y)]\right)\diff x \\
= \int_\R(g_1 + \I\:g_2)(t)(1+t^2)\diff\nu(t).
\end{multline*}
Observing that
\begin{multline*}
\Im[\qRe{q}(x+\I\:y)]+\I\:\Im[\qIm{q}(x+\I\:y)] = \frac{1}{2\I}(q(x+\I\:y)-\bar{\bar{q}(x+\I\:y)}) \\
= \frac{1}{2\I}(q(x+\I\:y)-q(x-\I\:y))
\end{multline*}
finishes the proof.
\endproof

While not directly related to the Stieltjes inversion formula \eqref{eq:quasi_stieltjes_v1}, the following result distils important additional information about the representing measure $\nu$ from Theorem \ref{thm:intRep_1var_quasi}.

\begin{prop}\label{prop:quasi_point_mass_limits}
Let $q$ be a quasi-Herglotz function, let $\nu$ be its representing measure in the sense of Theorem \ref{thm:intRep_1var_quasi} and let $m \in \N$. Then, for any $t_0 \in \R$, it holds that
$$\lim\limits_{z \ntto t_0}(t_0-z)^mq(z) = \lim\limits_{z \nntto t_0}(t_0-z)^mq(z) = \left\{\begin{array}{rcl}
(1+t_0^2)\nu(\{t_0\}) & ; & m = 1, \\
0 & ; & m \geq 2.
\end{array}\right.$$
\end{prop}

The proof of this results follows from the corresponding results for ordinary Herglotz functions, see \eg \cite{KacKrein1974}, and is, hence, omitted here.

\subsection{Zeros and compositions}

In some sense, quasi-Herglotz functions behave very similarly to ordinary Herglotz functions, whereas in other respects, they are quite different. For example, it is well known that Herglotz-functions have neither poles nor zeros in $\C\setminus\R$. This follows from their definition and the fact that for any Herglotz function $h$ also the quotient $z\mapsto-\frac1{h(z)}$ is a Herglotz-function.

For quasi-Herglotz functions, however, the situation is different. By definition, their poles are restricted to the real line, however, they can have non-real zeros of arbitrary order. Hence, for a quasi-Herglotz function $q$, in general, rational transformations of $q$ are {no longer} quasi-Herglotz-functions.

More generally speaking, in the case of ordinary Herglotz functions, it is an immediate consequence of the maximum principle that a Herglotz function $h$ attains a real value at some point in $\C\setminus\R$ if and only if the function $h$ is identically equal to a real-constant function. Hence, one may always compose two Herglotz functions as long as the first function is not a real-constant function. This is not true anymore for quasi-Herglotz functions in general. Consider, for example, the compositions
$$z \mapsto -\frac{1}{h(z)} \quad\text{and}\quad z \mapsto h\big(-\tfrac{1}{z}\big).$$
If $q$ is a quasi-Herglotz function, then the function
$${z \mapsto q\big(-\tfrac{1}{z}\big)}$$
is still a quasi-Herglotz function, which follows from the fact that every quasi-Herglotz function can be written in the form \eqref{eq:quasi_combination} together with the corresponding property for ordinary Herglotz functions. On the other hand, the function
$${z \mapsto - \tfrac{1}{q(z)}}$$
will not be well-defined as soon as the function $q$ attains a zero in $\C\setminus\R$.

\section{An analytic characterization of quasi-Herglotz functions}\label{sec:analytic}

The following theorem answers the question which holomorphic functions on $\C\setminus\R$ can be written as a (complex) linear combination of Herglotz functions by giving an analytic characterization of quasi-Herglotz functions.

\begin{thm}\label{thm:quasi_Tumarkin_Vladimirov}
Let $q \colon \C\setminus\R \to \C$ be a holomorphic function. Then $q$ is a quasi-Herglotz function if and only if the function $q$ satisfies, first, a growth condition, namely that there exists a number $M \geq 0$ such that
\begin{equation}
    \label{eq:Vladimirov_growth}
    |q(z)| \leq M \frac{1+|z|^2}{|\Im[z]|},
\end{equation}
for all $z \in \C\setminus\R$ and, second, the regularity condition
\begin{equation}\label{eq:Tumarkin_condition}
\sup_{r \in (0,1)}\int_\R|q(t+\I\:r)-q(t-\I\:r)|\frac{\diff t}{1+t^2} < \infty.
\end{equation}
\end{thm}

\proof
\emph{Step 1:} Assume first that the function $q$ satisfies conditions \eqref{eq:Vladimirov_growth} and \eqref{eq:Tumarkin_condition}. Following the idea of Vladimirov's proof of the integral representation theorem for Herglotz functions, \ie Theorem \ref{thm:intRep_1var}, \cite[pp. 290--292]{Vladimirov1979}, we are going to show that the function $q$, under our assumptions, admits a representation of the form \eqref{eq:intRep_1var_quasi} and is, hence, a quasi-Herglotz function by Theorem \ref{thm:intRep_1var_quasi}.

For any $r \in (0,1)$, consider the auxiliary functions $f_r$ and $g_r$ defined by
$$f_r(z) := \frac{q(z+\I\:r)}{1+z^2} \quad\text{and}\quad g_r(z) := \frac{q(z - \I\:r)}{1+z^2}.$$
The function $f_r$ is meromorphic on $\C\setminus\{\Im[z] = -r\}$, while the function $g_r$ is meromorphic on $\C\setminus\{\Im[z] = r\}$. Both functions have simple poles at $\pm\I$.

\begin{figure}[!ht]
\begin{tikzpicture}
\def\bigradius{2}
\def\shift{6.8}

\draw [help lines,->] (-1.25*\bigradius, 0) -- (1.25*\bigradius,0) node[above] {$x$};
\draw [help lines,->] (0, -1.25*\bigradius) -- (0, 1.25*\bigradius) node[right] {$\I\,y$};
\draw[line width=1pt,   decoration={ markings,
  mark=at position 0.1455 with {\arrow[line width=1.2pt]{>}},
  mark=at position 0.73 with {\arrow[line width=1.2pt]{>}}},
  postaction={decorate}]
  (\bigradius,0) arc (0:180:\bigradius) -- cycle;
\node at (-0.9*\bigradius,0.9*\bigradius) {$\Theta_R^+$};
\fill (0,\bigradius/2)  circle[radius=1.5pt];
\node at (0.25,\bigradius/2) {$+\I$};
\fill (0,-\bigradius/2)  circle[radius=1.5pt];
\node at (0.25,-\bigradius/2) {$-\I$};

\draw [help lines,->] (-1.25*\bigradius+\shift, 0) -- (1.25*\bigradius+\shift,0) node[above] {$x$};
\draw [help lines,->] (0+\shift, -1.25*\bigradius) -- (0+\shift, 1.25*\bigradius) node[right] {$\I\,y$};
\draw[line width=1pt,   decoration={ markings,
  mark=at position 0.4455 with {\arrow[line width=1.2pt]{<}},
  mark=at position 0.87 with {\arrow[line width=1.2pt]{<}}},
  postaction={decorate}]
  (-\bigradius+\shift,0) arc (180:360:\bigradius) -- cycle;
\node at (-0.9*\bigradius+\shift,-0.9*\bigradius) {$\Theta_R^-$};
\fill (0+\shift,\bigradius/2)  circle[radius=1.5pt];
\node at (0.25+\shift,\bigradius/2) {$+\I$};
\fill (0+\shift,-\bigradius/2)  circle[radius=1.5pt];
\node at (0.25+\shift,-\bigradius/2) {$-\I$};
\end{tikzpicture}
    
\caption{The contours of integration $\Theta_R^+$ (left) and $\Theta_R^-$ (right) with respect ot the points $\pm\I$.}
\label{fig:contours}
\end{figure}
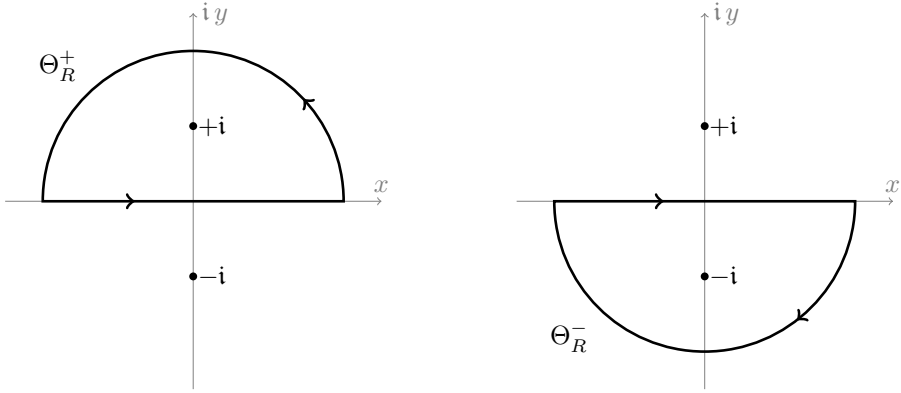

Let now $R > 1$ and take $\Theta_R^+$ to be the standard upper half-circle contour in $\C$, consisting of the line segment $[-R,R]$ and the half-circle $\theta_R^+$, and oriented counter-clockwise. Similarly, let $\Theta_R^-$ be the standard lower half-circle contour in $\C$, consisting of the line segment $[-R,R]$ and the half-circle $\theta_R^-$, and oriented clockwise. The contours $\Theta_R^+$ and $\Theta_R^-$ are visualized in Figure \ref{fig:contours}. 

\emph{Step 1.A:} Assume from now on that $z \in \C^+$ and that $R > |z|$. Then, by the residue theorem, it holds that
\begin{multline*}
    \frac{1}{2\pi\I}\oint_{\Theta_R^+}\frac{f_r(t)}{t-z}\diff t = \Res\left(\xi \mapsto \frac{f_r(\xi)}{\xi-z};z\right) + \Res\left(\xi \mapsto \frac{f_r(\xi)}{\xi-z};\I\right) \\
    = f_r(z) + \lim\limits_{\xi \to \I}\frac{q(\xi + \I\:r)}{(\xi-\I)(\xi + \I)(\xi - z)}(\xi - \I) = f_r(z) - \frac{q(\I + \I\:r)}{2\I(z-\I)}
\end{multline*}
and
\begin{multline*}
    \frac{1}{2\pi\I}\oint_{\Theta_R^-}\frac{g_r(t)}{t-z}\diff t = -\Res\left(\xi \mapsto \frac{g_r(\xi)}{\xi-z};-\I\right) \\
    = -\lim\limits_{\xi \to -\I}\frac{q(\xi - \I\:r)}{(\xi-\I)(\xi + \I)(\xi - z)}(\xi + \I) = -\frac{q(-\I - \I\:r)}{2\I(z + \I)}.
\end{multline*}

Using inequality \eqref{eq:Vladimirov_growth}, we may now estimate the integrals of the functions $f_r$ and $g_r$ over the arcs $\theta_\R^+$ and $\theta_R^-$, respectively. In particular, it holds that
\begin{multline*}
    \left|\int_{\theta_R^+}\frac{f_r(t)}{t-z}\diff t\right| = \left|\int_{\theta_R^+}\frac{q(t+\I\:r)}{(t-z)(1+t^2)}\diff t\right| \\
    = \left|\int_0^\pi\frac{q(R\E^{\I\theta}+\I\:r)\:\I\:R\E^{\I\theta}}{(R\E^{\I\theta}-z)(1+R^2\E^{2\I\theta})}\diff \theta\right| \leq \int_0^\pi \frac{R\:|q(R\E^{\I\theta}+\I\:r)|}{|R\E^{\I\theta}-z|\:|1+R^2\E^{2\I\theta}|}\diff\theta \\
    \leq  \frac{R\:M\:(1+(r+R)^2)}{(R-|z|)(R^2-1)}\int_0^\pi\frac{1}{R\sin\theta+r}\diff\theta,
\end{multline*}
where the last inequality holds due to the assumption that the function $q$ satisfies the growth condition \eqref{eq:Vladimirov_growth}. Using a standard trigonometric integral-substitution, we compute, for $R > r > 0$, that
\begin{multline*}
\int_0^\pi\frac{1}{R\sin\theta+r}\diff\theta = \left[\frac{1}{\sqrt{R^2-r^2}}\ln\left|\frac{\sqrt{R^2-r^2}-R-r \tan(\frac{\theta}{2})}{\sqrt{R^2-r^2}+R+r \tan(\frac{\theta}{2})}\right|\:\right]_{\theta = 0}^{\theta = \pi} \\
= \frac{1}{\sqrt{R^2-r^2}}\ln\left|\frac{\sqrt{R^2-r^2}-R}{\sqrt{R^2-r^2}+R}\right| \xrightarrow{R \to \infty} 0.
\end{multline*}
Therefore, we conclude that
\begin{equation}
    \label{eq:arc_estimate_f}
    \left|\int_{\theta_R^+}\frac{f_r(t)}{t-z}\diff t\right| \xrightarrow{R \to \infty} 0.
\end{equation}
Analogously, one may show that
\begin{equation}
    \label{eq:arc_estimate_g}
    \left|\int_{\theta_R^-}\frac{g_r(t)}{t-z}\diff t\right| \xrightarrow{R \to \infty} 0,
\end{equation}
allowing us to conclude that
\begin{equation}
    \label{eq:arc_estimate_combine}
    \lim\limits_{R \to \infty}\oint_{\Theta_R^+}\frac{f_r(t)}{t-z}\diff t = \int_\R\frac{f_r(t)}{t-z}\diff t \quad\text{and}\quad \lim\limits_{R \to \infty}\oint_{\Theta_R^-}\frac{g_r(t)}{t-z}\diff t = \int_\R\frac{g_r(t)}{t-z}\diff t.
\end{equation}

Combining these results with the previous calculations yields
\begin{eqnarray*}
f_r(z) & = & \frac{1}{2\pi\I}\int_\R\frac{f_r(t)}{t-z}\diff t + \frac{q(\I + \I\:r)}{2\I(z-\I)}, \\
0 & = & \frac{1}{2\pi\I}\int_\R\frac{g_r(t)}{t-z}\diff t + \frac{q(-\I - \I\:r)}{2\I(z + \I)}.
\end{eqnarray*}
Subtracting the second of the above equalities from the first and multiplying both sides of the result by $1+z^2$ yields
\begin{eqnarray}\label{eq:Vladimirov_intRep_temp}
q(z + \I\:r) & = & \frac{1+z^2}{2\pi\I}\int_\R\frac{1}{t-z}(f_r(t)-g_r(t))\diff t \\[0.25cm]
~ & ~ & + (1+z^2)\frac{q(\I + \I\:r)}{2\I(z-\I)} - (1+z^2)\frac{q(-\I - \I\:r)}{2\I(z + \I)} \nonumber \\[0.25cm]
~ & = & \frac{1+z^2}{2\pi\I}\int_\R\frac{1}{t-z}(q(t + \I\:r) - q(t - \I\:r))\frac{1}{1+t^2}\diff t \nonumber \\[0.25cm]
~ & ~ & + \frac{z+\I}{2\I}q(\I + \I\:r) - \frac{z- \I}{2\I}q(-\I - \I\:r) \nonumber.
\end{eqnarray}
Let $(r_n)_{n\in\N}\subseteq (0,1)$ be a monotonically decreasing sequence converging to zero. For any $n \in \N$, we define now a complex Borel measure $\nu_n$ on $\R$ via
$$\diff\nu_n(t) := \frac{1}{2\I}(q(t + \I\:r_n) - q(t - \I\:r_n))\frac{1}{1+t^2}\diff t.$$
By the assumption that the function $q$ satisfies the regularity condition \eqref{eq:Tumarkin_condition}, we infer that the sequence of measures $\{\nu_n\}_{n \in \N}$ is uniformly bounded in the natural norm on the space of complex Borel measures. Therefore, by Helly's selection principle \cite[Sec. VIII.4]{Natanson1964}, there exists a subsequence $\{r_{n_k}\}_{k \in \N} \subseteq (r_n)_{n \in \N}$ and a complex Borel measure $\nu$, such that
$$\nu_{n_k} \xrightarrow{w^*} \nu \quad\text{as}\quad k \to \infty,$$
where $w^*$ denotes that the limit is taken in the weak$^*$-sense.

As such, when taking the limit as $r \to 0^+$ in representation \eqref{eq:Vladimirov_intRep_temp}, we get
\begin{eqnarray}\label{eq:Vladimirov_intRep_temp_v3}
q(z) & = & \frac{1+z^2}{\pi}\int_\R\frac{1}{t-z}\diff\nu(t) + \frac{z+\I}{2\I}q(\I) - \frac{z- \I}{2\I}q(-\I)  \\[0.25cm]
~ & = & \frac{1+z^2}{\pi}\int_\R\frac{1}{t-z}\diff\nu(t) + \frac{z}{2\I}(q(\I)-q(-\I)) + \frac{1}{2}(q(\I)+q(-\I)). \nonumber
\end{eqnarray}
Furthermore, the integral in the above representation may be rewritten as
$$\frac{1+z^2}{\pi}\int_\R\frac{1}{t-z}\diff\nu(t) = \frac{1}{\pi}\int_\R\frac{1+t\:z}{t-z}\diff\nu(t) - \frac{z}{\pi}\int_\R\diff\nu(t),$$
refining representation \eqref{eq:Vladimirov_intRep_temp_v3} into
\begin{multline*}
q(z) = \frac{1}{\pi}\int_\R\frac{1+t\:z}{t-z}\diff\nu(t) + z\left(\frac{1}{2\I}(q(\I)-q(-\I)) - \frac{1}{\pi}\int_\R\diff\nu(t)\right) \\
+ \frac{1}{2}\big(q(\I)+q(-\I)\big).
\end{multline*}
This gives an integral representation of the form \eqref{eq:intRep_1var_quasi} for $z \in \C^+$.

\emph{Step 1.B:} On the other hand, if $z \in \C^-$ and $R > |z|$, it holds, by the residue theorem, that
\begin{multline*}
    \frac{1}{2\pi\I}\oint_{\Theta_R^+}\frac{f_r(t)}{t-z}\diff t = \Res\left(\xi \mapsto \frac{f_r(\xi)}{\xi-z};\I\right) \\
    = \lim\limits_{\xi \to \I}\frac{q(\xi + \I\:r)}{(\xi-\I)(\xi + \I)(\xi - z)}(\xi - \I) = - \frac{q(\I + \I\:r)}{2\I(z-\I)}
\end{multline*}
and
\begin{multline*}
    \frac{1}{2\pi\I}\oint_{\Theta_R^-}\frac{g_r(t)}{t-z}\diff t = -\Res\left(\xi \mapsto \frac{g_r(\xi)}{\xi-z};z\right) -\Res\left(\xi \mapsto \frac{g_r(\xi)}{\xi-z};-\I\right) \\
    = - g_r(z) -\lim\limits_{\xi \to -\I}\frac{q(\xi - \I\:r)}{(\xi-\I)(\xi + \I)(\xi - z)}(\xi + \I) = -g_r(z) - \frac{q(-\I - \I\:r)}{2\I(z + \I)}.
\end{multline*}

The estimates \eqref{eq:arc_estimate_f} and \eqref{eq:arc_estimate_g} hold also for $z \in \C^-$, which implies the equalities \eqref{eq:arc_estimate_combine}. Combining this result with the previous calculations yields
\begin{eqnarray*}
0 & = & \frac{1}{2\pi\I}\int_\R\frac{f_r(t)}{t-z}\diff t + \frac{q(\I + \I\:r)}{2\I(z-\I)}, \\
-g_r(z) & = & \frac{1}{2\pi\I}\int_\R\frac{g_r(t)}{t-z}\diff t + \frac{q(-\I - \I\:r)}{2\I(z + \I)}.
\end{eqnarray*}
Subtracting the second of the above equalities from the first and multiplying both sides of the result by $1+z^2$ yields
\begin{eqnarray*}\label{eq:Vladimirov_intRep_temp_v2}
q(z - \I\:r) & = & \frac{1+z^2}{2\pi\I}\int_\R\frac{1}{t-z}(q(t + \I\:r) - q(t - \I\:r))\frac{1}{1+t^2}\diff t \\[0.25cm]
~ & ~ & + \frac{z+\I}{2\I}q(\I + \I\:r) - \frac{z- \I}{2\I}q(-\I - \I\:r).
\end{eqnarray*}
Taking $r \to 0^+$ gives a representation of the form \eqref{eq:intRep_1var_quasi} for $z \in \C^-$.

\emph{Step 1.C:} The previous two steps have yielded, separately on $\C^+$ and $\C^-$, an integral representation of the function $q$. However, the parameters of both representations are identical, yielding a unified representation of the from \eqref{eq:intRep_1var_quasi} on $\C\setminus\R$. This finishes the first part of the proof.

\emph{Step 2:} Conversely, assume now that the function $q$ is a quasi-Herglotz function. Then there exists four ordinary Herglotz functions $h_1,h_2,h_3$ and $h_4$ such that equality \eqref{eq:quasi_combination} is satisfied.

By a result of Vladimirov, \cf \cite[Thm. 1]{Vladimirov1969} and \cite[pg. 203]{Vladimirov1979}, ordinary Herglotz functions satisfy the growth condition \eqref{eq:Vladimirov_growth} for $z \in \C^+$. The symmetry formula for ordinary Herglotz functions extends this result to $\C^-$. Indeed, if a Herglotz function $h$ satisfies condition \eqref{eq:Vladimirov_growth} for $z \in \C^+$ with some constant $M \geq 0$, then, for $z \in \C^-$, it holds that
$$|h(z)| = |h(\bar{z})| \leq M\frac{1+|z|^2}{-\Im[z]} = M\frac{1+|z|^2}{|\Im[z]|}.$$
Therefore, as any quasi-Herglotz function can be written in the form \eqref{eq:quasi_combination}, it will also satisfy condition \eqref{eq:Vladimirov_growth}.

In order to show that the function $q$ also satisfies the regularity condition \eqref{eq:Tumarkin_condition}, we calculate first using representation \eqref{eq:intRep_1var_quasi} that
$$q(z) - q(\bar{z}) = 2\,\I\:b\:\Im[z] + \frac{2\,\I}{\pi}\int_\R(1+t^2)\pois(z,t)\diff\nu(t)$$
for every $z \in \C\setminus\R$, where $\pois$ denotes the Poisson kernel as defined in formula \eqref{eq:poisson_kernels}. Hence, the integral appearing in condition \eqref{eq:Tumarkin_condition} may be estimated, for $r \in (0,1)$, as
\begin{eqnarray*}
0 & \leq & \int_\R|q(t+\I\:r)-q(t-\I\:r)|\frac{\diff t}{1+t^2} \\
~ & \leq & \int_\R\frac{2\,b\,r}{1+t^2}\diff t + \frac{2}{\pi}\int_\R\frac{1}{1+t^2}\left(\int_\R(1+\tau^2)\pois(t+\I\:r,\tau)\diff|\nu|(\tau)\right)\diff t \\
~ & \leq & 2\,b\,\pi + \frac{2}{\pi}\int_\R(1+\tau^2)\left(\int_\R\frac{\pois(t+\I\:r,\tau)}{1+t^2}\diff t\right)\diff|\nu|(\tau) \\
~ & = & 2\,b\,\pi + 2\int_\R\frac{(1+\tau^2)(1+r)}{(1+r)^2+\tau^2}\diff|\nu|(\tau).
\end{eqnarray*}
Here, the Fubini-Tonelli theorem was used to change the order of integration between the $t$- and $\tau$-variables. We observe now that
$$\frac{1}{2} \leq \frac{(1+\tau^2)(1+r)}{(1+r)^2+\tau^2} \leq 2$$
for all $\tau \in \R$ and all $r \in [0,1]$, allowing us to conclude that
$$\int_\R|q(t+\I\:r)-q(t-\I\:r)|\frac{\diff t}{1+t^2} \leq 2\,b\,\pi + 4\int_\R\diff|\nu|(\tau)$$
independently of $r \in (0,1)$, finishing the proof.
\endproof

\begin{remark}
In the regularity condition \eqref{eq:Tumarkin_condition}, it would suffice to assume the existence of a number $\delta > 0$ such that the supremum in condition \eqref{eq:Tumarkin_condition} is finite when taken over all $r \in (0, \delta)$.
\end{remark}

The preceeding theorem and its proof yield {additional information about the data of a quasi-Herglotz function} as well as the following corollaries.

\begin{coro}
Let $q$ be a quasi-Herglotz function given by the data $(a,b,\nu)$. Then, the number $b$ satisfies
$$b = \frac{1}{2\I}(q(\I)-q(-\I)) - \frac{1}{\pi}\int_\R\diff\nu(t)$$
and the measure $\nu$ is equal to the weak$^*$ limit, as $y \to 0^+$, of the measures given by the densities
$$x + \I\:y \mapsto \frac{1}{2\I}(q(x + \I\:y)-q(x - \I\:y))\frac{1}{1+x^2}.$$
\end{coro}

\begin{coro}
If $q \colon \C\setminus\R \to \C$ is a holomorphic function satisfying the growth condition \eqref{eq:Vladimirov_growth}, the regularity condition \eqref{eq:Tumarkin_condition} and, in addition, $q(\I) = q(-\I) = 0$, then there exists a complex Borel measure $\nu$ such that
$$q(z) = \frac{1+z^2}{\pi}\int_\R\frac{1}{t-z}\diff\nu(t).$$
\end{coro}

Furthermore, we may adapt Theorem \ref{thm:quasi_Tumarkin_Vladimirov} in order to give an analytic characterization of real quasi-Herglotz functions.

\begin{coro}
Let $q\colon \C\setminus\R \to \C$ be a holomorphic function. Then, $q$ is a real quasi-Herglotz function if and only if it satisfies the growth condition \eqref{eq:Vladimirov_growth}, the regularity condition \eqref{eq:Tumarkin_condition} and the symmetry condition
\begin{equation}
    \label{eq:temp_symm}
    q(z) = \bar{q(\bar{z})}
\end{equation}
for all $z \in \C\setminus\R.$
\end{coro}

\proof
If $q$ is a real quasi-Herglotz function, then it satisfies the three given conditions due to Theorems \ref{thm:intRep_1var_quasi} and \ref{thm:quasi_Tumarkin_Vladimirov}. Conversely, if the function satisfies the growth condition \eqref{eq:Vladimirov_growth} and the regularity condition \eqref{eq:Tumarkin_condition}, it is a quasi-Herglotz function by Theorem \ref{thm:quasi_Tumarkin_Vladimirov}. Therefore, it remains to show that the symmetry condition \eqref{eq:temp_symm} implies that all of the representing parameters of the function $q$ are real.

For the number $a$, we infer, using Corollary \ref{coro:a_and_b_parameters}, that
$$a = \frac{1}{2}\big(q(\I) + q(-\I)\big) = \frac{1}{2}\big(q(\I) + \bar{q(\I)}\big) = \Re[q(\I)] \in \R.$$
Similarly, for the number $b$, Corollary \ref{coro:a_and_b_parameters} implies that
$$\bar{b} = \lim\limits_{z \ntto \infty}\frac{\bar{q(z)}}{\bar{z}} = \lim\limits_{z \ntto \infty}\frac{q(\bar{z})}{\bar{z}} = \lim\limits_{z \nntto \infty}\frac{q(z)}{z} = b,$$
yielding that $b \in \R$. Finally, note that, in this case,
$$\frac{1}{2\,\I}\big(q(x+\I\,y) - q(x-\I\,y)\big) = \frac{1}{2\,\I}\big(q(x+\I\,y) - \bar{q(x+\I\,y)}\big) = \Im[q(x+\I\,y)],$$
yielding, via Proposition \ref{prop:Stieltjes_inversion}, that the measure $\nu$ is a signed measure. This finishes the proof.
\endproof

The growth condition \eqref{eq:Vladimirov_growth} and the regularity condition \eqref{eq:Tumarkin_condition} are independent, complementary, conditions, as will be illustrated shortly via three examples. One may interpret the regularity condition \eqref{eq:Tumarkin_condition} as guaranteeing that the function $q$ behaves sufficiently well at the real axis, \cf Example \ref{ex:Tum_point}, as well as at infinity when when approaching tangentially to the real line, \cf Example \ref{ex:Tum_inf}. On the other hand, the growth condition \eqref{eq:Vladimirov_growth} guarantees that the function $q$ behaves sufficiently well at infinity when approaching non-tangentially, \cf Example \ref{ex:Vlad_inf}. In total, the two conditions constrain the behaviour of the function $q$ over all of the boundary of $\C\setminus\R$.

\begin{example}\label{ex:Tum_inf}
Consider the function
$$f(z) := \left\{\begin{array}{rcl}
z & ; & z \in \C^+, \\
-z & ; & z \in \C^-.
\end{array}\right.$$
This function satisfies condition \eqref{eq:Vladimirov_growth}, but does not satisfy condition \eqref{eq:Tumarkin_condition}. Indeed, we calculate that
$$|f(t+\I\:r)-f(t-\I\:r)|\:\frac{1}{1+t^2} = \frac{2|t|}{1+t^2},$$
which is not integrable at $\pm\infty$. \hfill$\lozenge$
\end{example}

\begin{example}\label{ex:Tum_point}
Consider the function
$$f(z) := \left\{\begin{array}{rcl}
-\frac{1}{z} & ; & z \in \C^+, \\[0.2cm]
\frac{1}{z} & ; & z \in \C^-.
\end{array}\right.$$
Also this function satisfies condition \eqref{eq:Vladimirov_growth}, but not condition \eqref{eq:Tumarkin_condition}. Indeed, we calculate that
$$\int_\R|f(t+\I\:r)-f(t-\I\:r)|\:\frac{1}{1+t^2}\diff t = 4\int_{0}^\infty\frac{t}{(1+t^2)(r^2+t^2)}\diff t = \frac{4\ln r}{r^2-1},$$
which blows up as $r \to 0^+$. However, contrary to the previous example, the integrability problem does not lie at $\pm \infty$ but at zero. This can also be seen by noting that for any sequence $(r_n)_{n \in \N} \subseteq (0,1)$ that monotonically decreases to zero the sequence of functions
$$|f(t+\I\:r_n)-f(t-\I\:r_n)|\:\frac{1}{1+t^2} = \frac{|2\,t|}{(1+t^2)(r_n^2+t^2)}$$
monotonically increases point-wise on $\R\setminus\{0\}$ to
$$\frac{2}{|t|(1+t^2)},$$
which is not integrable at $0$.
\hfill$\lozenge$
\end{example}

\begin{example}\label{ex:Vlad_inf}
Consider the function $f(z) := \E^{-z^2}$ for $z \in \C\setminus\R$. This function satisfies condition \eqref{eq:Tumarkin_condition}, but does not satisfy condition \eqref{eq:Vladimirov_growth}. Indeed, we calculate that
$$|f(t+\I\:r)-f(t-\I\:r)|\:\frac{1}{1+t^2} = 2\,\E^{-t^2+r^2}|\sin(2\,r\,t)|\frac{1}{1+t^2} \leq 2\,\E^{-t^2+1}\frac{1}{1+t^2}$$
which is integrable. On the other hand, for $z = \I\:y$, condition \eqref{eq:Vladimirov_growth} is equivalent to the inequality
$$|y|\E^{y^2} \leq M(1+|y|^2).$$
for some $M \in \R$. However, such an inequality can never be satisfied for large $y$.\hfill$\lozenge$
\end{example}

\section{Symmetry and uniqueness}\label{subsec:sym_and_uni}

We return now to the definition of Herglotz functions. As mentioned in Section \ref{sec:background}, Definition \ref{def:Herglotz} is not the only way Herglotz functions are defined in the literature. One way is to consider functions that are defined, \emph{a priori}, only on the upper half-plane $\C^+$.

\begin{define}\label{def:Herglotz_v2}
A function $h \colon \C^+ \to \C$ is called a \emph{Herglotz function (on $\C^+$)} if it is holomorphic with $\Im[h(z)] \geq 0$ for all $z \in \C^+$.
\end{define}

Using this definition, one may establish an integral representation as in Theorem \ref{thm:intRep_1var} for $z \in \C^+$, including the statement about uniqueness of the data. However, representation \eqref{eq:intRep_1var} is automatically well-defined for any $z \in \C\setminus\R$, which may be used to extend any Herglotz function on $\C^+$ to a function defined on $\C\setminus\R$. This extension is called the \emph{symmetric extension}, since it automatically satisfies {the condition that} $h(\overline z)=\overline{h(z)}$ and will be a Herglotz function on $\C\setminus\R$, \ie it will satisfy Definition \ref{def:Herglotz}. Conversely, the restriction to $\C^+$ of any Herglotz function on $\C\setminus\R$ will satisfy Definition \ref{def:Herglotz_v2}.

Hence a Herglotz-function is uniquely determined by its values in one half-plane only. However, as we will see, this is not true for quasi-Herglotz functions.

\begin{coro}\label{coro:quasi_symmetry}
Let $q$ be a quasi-Herglotz function and $\bar{q}$ its conjugate function. Then, it holds for every $z \in \C\setminus\R$ that
\begin{equation}
    \label{eq:quasi_symmetry}
    \bar{q(z)} = \bar{q}(\bar{z}).
\end{equation}
\end{coro}

\proof
Using, representation \eqref{eq:intRep_1var_quasi}, we calculate, for every $z \in \C\setminus\R$, that
$$\bar{q(z)} = \bar{a+b\:z + \frac{1}{\pi}\int_\R \til{K}(z,t)\diff\nu(t)} = \bar{a}+\bar{b}\:\bar{z} + \frac{1}{\pi}\int_\R \til{K}(\bar{z},t)\diff\bar{\nu}(t) = \bar{q}(\bar{z}),$$
as desired.
\endproof

For real quasi-Herglotz functions, Corollary \ref{coro:quasi_symmetry} implies immediately that values of the function in one half-plane are uniquely determined by its values in the other half-plane. Indeed, as $\bar{q} = q$ for real quasi-Herglotz functions, one may make use of the same idea as with ordinary Herglotz functions.

However, for general quasi-Herglotz functions, this is not necessarily the case. For example, it is not enough to say that we are considering the function $q(z) := \I$ for $z \in \C^+$, as it is unclear to which quasi-Herglotz function we are referring. It may be the quasi-Herglotz function $q(z) = \I$ for $z \in \C\setminus\R$, which is represented by the data $(\I,0,0)$ in the sense of Theorem \ref{thm:intRep_1var_quasi}, or it may be the function 
\begin{equation}
    \label{eq:Lebesgue_function}
    q(z) := \left\{\begin{array}{rcl}
    \I & ; & z \in \C^+, \\
    -\I & ; & z \in \C^-,
    \end{array}\right.
\end{equation}
which is represented by the data $(0,0,\til{\lambda}_\R)$, where $\til{\lambda}_\R$ denotes the finite Lebesgue measure on $\R$, \ie the measure given by formula \eqref{eq:measure_make_finite} if $\mu = \lambda_\R$.

Combining these two functions, one may construct the function
$$q_1(z) := \left\{\begin{array}{rcl}
0 & ; & z \in \C^+, \\
-\I & ; & z \in \C^-,
\end{array}\right.$$
which is represented by the data $(-\frac{1}{2}\,\I,0,\frac{1}{2}\til{\lambda}_\R)$ in the sense of Theorem \ref{thm:intRep_1var_quasi}. Hence, if one adds (any constant multiple of) the function $q_1$ to a given quasi-Herglotz function, the values of this function in the upper half-plane are not going to change, while the values in the lower half-plane will.

\section{Quasi-Herglotz functions that are identically zero in one half-plane}\label{sec:identically_zero}

Consider now the following question: if we are given a quasi-Herglotz function, how many other quasi-Herglotz functions attain the same values in one half-plane while attaining different values in the other half-plane. Answering this question amounts to characterizing all quasi-Herglotz functions that are identically zero in one half-plane, which we will do now.

\subsection{Characterizations}

We begin by presenting a corollary of the symmetry formula \eqref{eq:quasi_symmetry} for quasi-Herglotz functions that are identically zero in the lower half-plane.

\begin{coro}\label{coro:conjugate_function_zero}
Let $q$ be a quasi-Herglotz function that is given as
$$q(z) := \left\{\begin{array}{rcl}
f(z) & ; & z \in \C^+, \\[0.2cm]
0 & ; & z \in \C^-,
\end{array}\right.$$
for some function $f\colon\C^+ \to \C$. Then, its conjugate function $\bar{q}$ is equal to
$$\bar{q}(z) = \left\{\begin{array}{rcl}
0 & ; & z \in \C^+, \\[0.2cm]
\bar{f(\bar{z})} & ; & z \in \C^-.
\end{array}\right.$$
\end{coro}

Quasi-Herg\-lotz functions that are identically zero in the lower half-plane may also be characterized in terms of their quasi-real and quasi-imaginary parts.

\begin{coro}\label{coro:tautological}
Let $q$ be a quasi-Herglotz function. Then, the function $q$ is identically zero in the lower half-plane if and only if its quasi-real and quasi-imaginary parts satisfy the relation
\begin{equation}\label{eq:quasi_real_imagnary_zero}
\qRe{q}(z) = -\I\,\qIm{q}(z)
\end{equation}
for $z \in \C^-$.
\end{coro}

The following theorem gives, instead, a characterization of quasi-Herglotz functions that are identically zero in the lower half-plane in terms of the data of the function in question.

\begin{thm}\label{thm:tautological}
A quasi-Herglotz function $q$ is identically zero in the lower half-plane $\C^-$ if and only if it holds, for its data $(a,b,\nu)$, that $b = 0$,
\begin{equation}\label{eq:a_tautological}
a = \frac{\I}{\pi}\int_\R\diff\nu(t)
\end{equation}
and
\begin{equation}\label{eq:nu_tautological}
\int_\R\frac{t-\I}{t-z}\diff\nu(t) = 0
\end{equation}
for all $z \in \C^-$.
\end{thm}

\begin{remark}
In other words, there exists a bijection between the sub-space of quasi-Herglotz functions that are identically zero in the lower half-plane and the space of complex Borel measures satisfying  condition \eqref{eq:nu_tautological}.
\end{remark}

\proof
First, assume that $q$ is identically zero in $\C^-$. Then, by formula \eqref{eq:b_constant}, it holds that
$$b = \lim\limits_{z \nntto \infty}\frac{q(z)}{z} = \lim\limits_{z \nntto \infty}0 = 0.$$
Furthermore, by formula \eqref{eq:a_constant}, it holds that
$$a = \frac{1}{2}\big(q(\I)+\underbrace{q(-\I)}_{=\,0}\big) = \frac{1}{2}q(\I) = \frac{1}{2}\left(a + \frac{\I}{\pi}\int_\R\diff\nu(t)\right),$$
yielding the desired description of the constant $a$. Finally, for any $z \in \C^-$, it holds, using representation \eqref{eq:intRep_1var_quasi}, that
$$0 = q(z) = \frac{\I}{\pi}\int_\R\diff\nu(t) + \frac{1}{\pi}\int_\R \til{K}(z,t)\diff\nu(t) = \frac{z+\I}{\pi}\int_\R\frac{t-\I}{t-z}\diff\nu(t),$$
yielding that
$$\int_\R\frac{t-\I}{t-z}\diff\nu(t) = 0$$
for all $z \in \C^-\setminus\{-\I\}$. However, since the function
$$z \mapsto \int_\R\frac{t-\I}{t-z}\diff\nu(t)$$
is holomorphic in $\C^-$, it must, by the identity principle, also equal zero when $z = -\I$.

Conversely, assume that the data of the function $q$ satisfies the prescribed conditions. Then, representation \eqref{eq:intRep_1var_quasi} can be rewritten as
$$q(z) = \frac{\I}{\pi}\int_\R\diff\nu(t) + \frac{1}{\pi}\int_\R \til{K}(z,t)\diff\nu(t) = \frac{z+\I}{\pi}\int_\R\frac{t-\I}{t-z}\diff\nu(t),$$
which is equal to zero for any $z \in \C^-$ by the condition on the measure $\nu$. This finishes the proof.
\endproof

\begin{remark}
From the above proof, it follows immediately that for quasi-Herglotz functions that are identically zero in $\C^+$, the conditions on the data $(a,b,\nu)$ become $b = 0$,
$$a = -\frac{\I}{\pi}\int_\R\diff\nu(t)$$
and
$$\int_\R\frac{t-\I}{t-z}\diff\nu(t) = 0$$
for all $z \in \C^+$.
\end{remark}

\begin{example}
Consider the complex Borel measure $\nu$ on $\R$, defined by
$$\nu(U) := \int_\R\frac{\chi_U(t)}{(t+\I)^2}\diff t,$$
where $U \subseteq \R$ is a Borel measurable set and $\chi_U$ denotes its characteristic function. It follows now by standard residue calculus that
$$\int_\R\frac{t-\I}{t-z}\diff\nu(t) = \int_\R\frac{t-\I}{t-z}\,\frac{\diff t}{(t+\I)^2} = 0$$
for every $z \in \C^-$. Furthermore, it holds that $\nu(\R) = 0$. Therefore, the data $(0,0,\nu)$ satisfies the assumptions of Theorem \ref{thm:tautological} and defines, via representation \eqref{eq:intRep_1var_quasi}, a quasi-Herglotz function that is identically zero in $\C^-$. This function $q$ is equal to
$$q(z) = \frac{1}{\pi}\int_\R\frac{1+t\,z}{t-z}\,\frac{\diff t}{(t+\I)^2} =  \left\{\begin{array}{RCL}
2\,\I+\frac{4}{z+\I} & ; & z \in \C^+, \\[0.3cm]
0 & ; & z \in \C^-,
\end{array}\right.$$
while its conjugate function $\bar{q}$ may be obtained directly from Corollary \ref{coro:conjugate_function_zero}.\hfill$\lozenge$
\end{example}

\subsection{Refinements and other properties}

The first result gives a necessary condition on the measure $\nu$ to be a representing measure of a quasi-Herglotz function vanishing in one half-plane.

\begin{prop}\label{prop:quasi_point_mass}
Let $q$ be a quasi-Herglotz functions that is identically zero in (at least) one half-plane and let $\nu$ be its representing measure in the sense of Theorem \ref{thm:intRep_1var_quasi}. Then, the measure $\nu$ cannot have any point-masses, \ie it holds that $\nu(\{t_0\}) = 0$ for all points $t_0 \in \R$.
\end{prop}

\proof
By Proposition \ref{prop:quasi_point_mass_limits}, it holds, for any point $t_0 \in \R$, that
$$(1+t_0^2)\nu(\{t_0\}) = \lim\limits_{z \ntto t_0}(t_0-z)q(z) = \lim\limits_{z \nntto t_0}(t_0-z)q(z).$$
As the function $q$ is identically zero in (at least) one half-plane, (at least) one of the above limits, and thus both, is equal to zero. This gives the desired result.
\endproof

The following corollary describes in more detail the role of signed measures as representing measures of quasi-Herglotz functions that are identically zero in one half-plane.

\begin{coro}\label{coro:signed_zero}
Let $q$ be a quasi-Herglotz function represented by the data $(a,0,\nu)$ that is identically zero in (at least) one half-plane. Then, the following statements hold.
\begin{itemize}
    \item[(a)]{It holds that $a = 0$ if and only if $\nu(\R) = 0$.}
    \item[(b)]{If $a = 0$, then $\qIm{\nu} \not\equiv 0$ unless $\nu \equiv 0$, \ie $\nu$ cannot be a signed measure unless it is identically zero.}
    \item[(c)]{If $\qIm{\nu}\equiv 0$, then $\Re[a] = 0$ and $\nu = \Im[a]\,\til{\lambda}_\R$.}
\end{itemize}
\end{coro}

\proof
Statement (a) is an obvious consequence of Theorem \ref{thm:tautological}. To prove statement (b), assume, without loss of generality, that the function $q$ is identically zero in the lower half-plane but not identically zero in the upper half-plane. If $q$ is represented by the data $(0,0,\nu)$, where $\nu$ is a signed measure on $\R$, then it is a real quasi-Herglotz function, implying that its values in $\C\setminus\R$ are determined uniquely by its values in one half-plane, \cf Section \ref{subsec:sym_and_uni}. Therefore, the function $q$ is either identically zero overall in $\C\setminus\R$ or not identically zero in either half-plane. Both cases lead to a contradiction.

Finally, to prove statement (c), assume that $\qIm{\nu}\equiv 0$. By Theorem \ref{thm:tautological}, we infer that the quasi-real and quasi-imaginary parts of the function $q$ are given by the data
$$\qRe{q} \sim (\Re[a],0,\qRe{\nu}) \quad\text{and}\quad \qIm{q} \sim (\Im[a],0,0),$$
respectively. For $z \in \C^-$, we may now deduce, via Corollary \ref{coro:tautological}, that
$$\qRe{q}(z) = -\I\,\qIm{q}(z) = -\I\,\Im[a].$$
Therefore, for $z \in \C^+$, it holds that
$$\qRe{q}(z) = \bar{\qRe{q}(\bar{z})} = \I\,\Im[a].$$
The desired result now follows by the uniqueness-statement of Theorem \ref{thm:intRep_1var_quasi}.
\endproof

\section{Rational quasi-Herglotz functions}\label{sec:rational}

A large class of rational functions are actually quasi-Herglotz functions. In this section, we are studying this class and, specifically, are going to show that any rational quasi-Herglotz function may be decomposed into a sum of three quasi-Herglotz functions of a very particular type, \cf Theorem \ref{thm:quasi_rational_3}. We start with an easy but useful observation.

\begin{remark}\label{rem:rational}
If a quasi-Herglotz function $q$ is rational in a half-plane, \ie it can be written, for $z$ in that half-plane, as $q(z)=\frac{P(z)}{Q(z)}$ for two coprime complex polynomials $P$ and $Q$, then the polynomials $P$ and $Q$ have to satisfy the following properties.
\begin{enumerate}
    \item[(i)]{It holds that $\deg(P) \leq \deg(Q)+1$ (due to the existence of the limit in equality \eqref{eq:b_constant}),}
    \item[(ii)]{The zeros of the polynomial $Q$ do not lie in this half plane (due to holomorphy).}
    \item[(iii)]{The real zeros of the polynomial $Q$ are simple (due to Proposition \ref{prop:quasi_point_mass_limits}).}
    \item[(iv)]{A point $t_0 \in \R$ is a zero of the polynomial $Q$ if and only if $\nu(\{t_0\}) \neq 0$ where $\nu$ denotes the representing measure of the function $q$ (also due to Proposition \ref{prop:quasi_point_mass_limits}).}
\end{enumerate}
\end{remark}

The first theorem of this section characterizes quasi-Herglotz functions that are equal to a rational functions in $\C^+$ and identically zero in $\C^-$.

\begin{thm}\label{thm:quasi_rational_1}
Let $q\colon \C\setminus\R \to \C$ be a holomorphic function for which there exist two coprime complex polynomials such that 
\begin{equation}
\label{eq:quasi_ratioanl_1}
q(z) = \left\{\begin{array}{RCL}
\frac{P(z)}{Q(z)} & ; & z \in \C^+, \\[0.3cm]
0 & ; & z \in \C^-.
\end{array}\right.
\end{equation}
Then, the function $q$ is a quasi-Herglotz function if and only if the polynomials $P$ and $Q$ are such that $Q(z) \neq 0$ for all $z \in \C^+\cup\R$ and $\deg(P) \leq \deg(Q)$.
\end{thm}

\proof
Assume first that the function $q$ of the form \eqref{eq:quasi_ratioanl_1} is a quasi-Herglotz function. By Corollary \ref{coro:a_and_b_parameters}, it holds that
$$0 = \lim\limits_{z \nntto \infty}\frac{q(z)}{z} = \lim\limits_{z \ntto \infty}\frac{q(z)}{z} = \lim\limits_{z \ntto \infty}\frac{P(z)}{z\,Q(z)},$$
which is only possible if $\deg(P) \leq \deg(Q)$.

Furthermore, by Remark \ref{rem:rational}.(ii) and \ref{rem:rational}.(iii), the polynomial $Q$ has, at most, simple real zeros. However, even this may not occur due to Proposition \ref{prop:quasi_point_mass} and \ref{rem:rational}.(iv), as desired.

Conversely, suppose that the polynomials $P$ and $Q$ satisfy the conditions of the theorem. Then, the function $q$ of the form \eqref{eq:quasi_ratioanl_1} clearly satisfies the growth condition \eqref{eq:Vladimirov_growth} due to the degree-constraint for the polynomials $P$ and $Q$. Therefore, it remains to show that the function $q$ also satisfies the regularity condition \eqref{eq:Tumarkin_condition}, as the result, thereafter, follows from Theorem \ref{thm:quasi_Tumarkin_Vladimirov}.

To that end, we note that, due to the assumptions on the polynomials $P$ and $Q$, we have, for any fixed $r \in [0,1]$, that the function
$$t \mapsto \left|\frac{P(t+ \I\,r)}{Q(t+\I\,r)}\right|$$
is a bounded continuous function on $\R$, with an upper bound that, in general, depends on $r \in [0,1]$. However, as the interval $[0,1]$ is compact, there exists an upper bound that is independent of $r$, \ie there exists a constant $C \in \R$ such that
$$\left|\frac{P(t+ \I\,r)}{Q(t+\I\,r)}\right| \leq C$$
for all $t \in \R$ and all $r \in [0,1]$. Hence,
$$\int_\R\left|\frac{P(t+ \I\,r)}{Q(t+\I\,r)}\right|\frac{1}{1+t^2}\diff t \leq C\,\pi,$$
finishing the proof.
\endproof

\begin{remark}\label{rem:quasi_ratioanl_1}
A trivial - but later needed - reformulation of the above theorem states that a holomorphic function $q\colon \C\setminus\R \to \C$ for which there exist two coprime complex polynomials such that 
$$q(z) = \left\{\begin{array}{RCL}
0 & ; & z \in \C^+, \\
\frac{P(z)}{Q(z)} & ; & z \in \C^-,
\end{array}\right.$$
is a quasi-Herglotz function if and only if the polynomials $P$ and $Q$ are such that $Q(z) \neq 0$ for all $z \in \C^-\cup\R$ and $\deg(P) \leq \deg(Q)$.
\end{remark}

The second theorem of this section describes which quasi-Herglotz functions are equal to rational functions on $\C$.

\begin{thm}\label{thm:quasi_rational_2}
Let $q\colon \C\setminus\R \to \C$ be a holomorphic function for which there exist two coprime complex polynomials such that 
\begin{equation}
\label{eq:quasi_ratioanl_2}
q(z) = \frac{P(z)}{Q(z)}
\end{equation}
for all $z \in \C\setminus\R$. Then, the function $q$ is a quasi-Herglotz function if and only if the polynomials $P$ and $Q$ are such that all the zeros of $Q$ are simple real zeros and $\deg(P) \leq \deg(Q) + 1$.
\end{thm}

\proof
Assume first that the function $q$ of the form \eqref{eq:quasi_ratioanl_2} is a quasi-Herglotz function. Then, by Remark \ref{rem:rational}, the polynomials $P$ and $Q$ satisfy all of the conditions listed in the theorem.

Conversely, suppose that the polynomials $P$ and $Q$ satisfy the conditions of the theorem. Due to the degree-constraint for the polynomials $P$ and $Q$, any function $q$ of the form \eqref{eq:quasi_ratioanl_2} clearly satisfies the growth condition \eqref{eq:Vladimirov_growth}. In order to show that any such function $q$ also satisfies the regularity condition \eqref{eq:Tumarkin_condition}, we note that, due to the assumptions on the polynomials $P$ and $Q$, there exist a number $b \in \C$ and a complex polynomial $\til{P}$ with $\deg(\til{P}) \leq \deg(Q)$ such that
$$\frac{P(z)}{Q(z)} = b\,z + \frac{\til{P}(z)}{Q(z)}$$
for all $z \in \C\setminus\R$. Thus, we have, for any fixed $r \in [0,1]$, that the function
\begin{multline*}
t \mapsto \left|b(t+\I\,r)+\frac{\til{P}(t+ \I\,r)}{Q(t+\I\,r)} - b(t-\I\,r)-\frac{\til{P}(t-\I\,r)}{Q(t-\I\,r)}\right| \\
= \left|2\,\I\,r\,b + \frac{\til{P}(t+\I\,r)Q(t-\I\,r)-\til{P}(t-\I\,r)Q(t+\I\,r)}{Q(t+\I\,r)Q(t-\I\,r)} \right|
\end{multline*}
is a bounded continuous function on $\R$, with an upper bound that, in general, depends on $r \in [0,1]$. The results now follows via an analogous argument as in the proof of the previous theorem.
\endproof

The third and final theorem of this section gives the announced decomposition of a general reational quasi-Herglotz function.

\begin{thm}\label{thm:quasi_rational_3}
Let $q\colon \C\setminus\R \to \C$ be a holomorphic function for which there exist two pairs of coprime complex polynomials $P_1,Q_1$ and $P_2,Q_2$ such that 
\begin{equation}\label{eq:quasi_rational_3}
q(z) := \left\{\begin{array}{RCL}
\frac{P_1(z)}{Q_1(z)} & ; & z \in \C^+, \\[0.3cm]
\frac{P_2(z)}{Q_2(z)} & ; & z \in \C^-.
\end{array}\right.
\end{equation}
Then, the function $q$ is a quasi-Herglotz function if and only if it can be written as a sum of quasi-Herglotz functions from Theorem \ref{thm:quasi_rational_1}, Remark \ref{rem:quasi_ratioanl_1} and Theorem \ref{thm:quasi_rational_2}.
\end{thm}

\proof
If a function can be written as a sum of functions from Theorem \ref{thm:quasi_rational_1}, Remark \ref{rem:quasi_ratioanl_1} and Theorem \ref{thm:quasi_rational_2} it is obviously a quasi-Herglotz function.

Conversely, assume that a function $q$ of the form \eqref{eq:quasi_rational_3} is a quasi-Herglotz function. By Corollary \ref{coro:a_and_b_parameters}, it holds that
$$\lim\limits_{z \ntto \infty}\frac{q(z)}{z} = \lim\limits_{z \ntto \infty}\frac{P_1(z)}{z\,Q_1(z)} < \infty$$
and
$$\lim\limits_{z \nntto \infty}\frac{q(z)}{z} = \lim\limits_{z \nntto \infty}\frac{P_2(z)}{z\,Q_2(z)} < \infty.$$
Furthermore, the two limits are always equal, implying that there exist numbers $a_1,a_2,b \in \C$ and complex polynomials $\til{P}_1,\til{Q}_1,\til{P}_2$ and $\til{Q}_2$ such that 
$$q(z) = b\,z + \left\{\begin{array}{RCL}
a_1 + \frac{\til{P}_1(z)}{\til{Q}_1(z)} & ; & z \in \C^+, \\[0.4cm]
a_2 + \frac{\til{P}_2(z)}{\til{Q}_2(z)} & ; & z \in \C^-,
\end{array}\right.$$
where $\til{P}_1$ and $\til{Q}_1$ are coprime with $\deg(\til{P}_1) < \deg(\til{Q}_1)$ and $\til{P}_2$ and $\til{Q}_2$ are coprime with $\deg(\til{P}_2) < \deg(\til{Q}_2)$.

Via the fundamental theorem of algebra, we may, for $j=1,2$, factorize the polynomial $\til{Q}_j$ as
$$\til{Q}_j(z) = Q_{j1}(z)Q_{j2}(z),$$
where the polynomial $Q_{j1}$ has leading coefficient one and only has real zeros, while the polynomial $Q_{j2}$ only has zeros lying in $\C^-$ when $j = 1$ and only has zeros lying in $\C^+$ when $j = 2$. Using partial fraction decompositions, the function $q$ may now be rewritten as
\begin{equation}
\label{eq:quasi_rational_3_decomp}
q(z) = b\,z + \left\{\begin{array}{R}
\frac{P_{11}(z)}{Q_{11}(z)} \\[0.4cm]
\frac{P_{21}(z)}{Q_{21}(z)} 
\end{array}\right. + \left\{\begin{array}{R}
a_1 + \frac{P_{12}(z)}{Q_{12}(z)} \\[0.8cm]
0
\end{array}\right. + \left\{\begin{array}{RCL}
0 & ; & z \in \C^+, \\[0.4cm]
a_2 + \frac{P_{22}(z)}{Q_{22}(z)} & ; & z \in \C^-,
\end{array}\right.
\end{equation}
where, for $j=1,2$, we have $P_{j1}(z)Q_{j2}(z)+P_{j2}(z)Q_{j1}(z) = \til{P}_j(z)$ with $\deg(P_{j1}) < \deg(Q_{j1})$ and $\deg(P_{j2}) < \deg(Q_{j2})$. Note also that, for $j = 1,2$, the pairs of polynomials $P_{j1}$ and $Q_{j1}$, as well as $P_{j2}$ and $Q_{j2}$, are coprime.

The two functions in the decomposition \eqref{eq:quasi_rational_3_decomp} that are identically zero in one half-plane clearly satisfy the conditions of Theorem \ref{thm:quasi_rational_1} or Remark \ref{rem:quasi_ratioanl_1}, respectively. Therefore, it remains to show that the function
$$q_1\colon z \mapsto \left\{\begin{array}{RCL}
\frac{P_{11}(z)}{Q_{11}(z)} & ; & z \in \C^+, \\[0.4cm]
\frac{P_{21}(z)}{Q_{21}(z)} & ; & z \in \C^-, 
\end{array}\right.  $$
satisfies the conditions of Theorem \ref{thm:quasi_rational_2}. To that end, note that the function $q_1$ is, for certain, a quasi-Herglotz function, as all of the other functions in decomposition \eqref{eq:quasi_rational_3_decomp} have been shown to, or are assumed to, be quasi-Herglotz functions.

Hence, we may apply Proposition \ref{prop:quasi_point_mass_limits}, yielding, for $m \geq 2$ and $t_0 \in \R$, that
$$0 = \lim\limits_{z \ntto t_0}(t_0 - z)^m\frac{P_{11}(z)}{Q_{11}(z)} = \lim\limits_{z \nntto t_0}(t_0 - z)^m\frac{P_{21}(z)}{Q_{21}(z)}.$$
Thus, both polynomials $Q_{11}$ and $Q_{21}$ only have simple zeros. For $m = 1$, the above limits are still equal, though they may be non-zero, with this option occurring if and only if one, and thus both, of the polynomials has a zero at $t_0$. Therefore, we conclude that the polynomials $Q_{11}$ and $Q_{12}$ have identical zeros while having the same leading coefficient. Therefore, they are the same.

As such, it remains to show that the polynomials $P_{11}$ and $P_{12}$ also are identical. To that end, denote $k = \deg(Q_{11})$ and let $t_1,t_2,\ldots,t_k$ be the zeros of $Q_{11}$, \ie
$$Q_{11}(z) = \prod_{\ell = 1}^k(z - t_\ell).$$
At any point $t_\ell$ with $1 \leq \ell \leq k$ it holds, due to Proposition \ref{prop:quasi_point_mass_limits}, that
$$P_{11}(t_\ell)\prod_{\substack{j = 1 \\ j \neq \ell}}^k(z - t_\ell)^{-1} = P_{12}(t_\ell)\prod_{\substack{j = 1 \\ j \neq \ell}}^k(z - t_\ell)^{-1}.$$
Hence, the polynomials $P_{11}$ and $P_{12}$ coincide in $k$ points. As their degrees are, in addition, strictly less than than $k$, the result follows.
\endproof

\begin{remark}
If one assumes, as in the proof above, that all of the $P,Q$-pairs of polynomials appearing in formula \eqref{eq:quasi_rational_3_decomp} are coprime, then the decomposition in formula \eqref{eq:quasi_rational_3_decomp} is unique for a given function $q$.
\end{remark}

Theorem \ref{thm:quasi_rational_3} implies that there exists no rational quasi-Herglotz function of the form
$$q(z) =\left\{\begin{array}{RCL}
\frac{P_{1}(z)}{Q(z)} & ; & z \in \C^+, \\[0.4cm]
\frac{P_{2}(z)}{Q(z)} & ; & z \in \C^-,
\end{array}\right.$$
with distinct polynomials $P_1$ and $P_2$ both having {degree grater or eqaul to 1}. Furthermore, for any rational quasi-Herglotz function written in the form \eqref{eq:quasi_rational_3}, the polynomials $Q_1$ and $Q_2$ must have equal real zeros, \ie there exists no quasi-Herglotz function for which $Q_1(t_0) = 0$ for some $t_0 \in \R$, but $Q_2(t_0) \neq 0$. Furthermore, as noted in Remark \ref{rem:rational}.(iii), all of the real poles of a rational quasi-Herglotz function must be simple. On the other hand, Theorem \ref{eq:intRep_1var_quasi} shows that quasi-Herglotz functions are the largest class of holomoprhic functions on $\C\setminus\R$ admitting a representation of the form \eqref{eq:intRep_1var_quasi}. Hence, we may ask what would need to change in representation \eqref{eq:intRep_1var_quasi} in order to allow for \eg rational functions with second order poles on $\R$. As the problem with higher order real poles stems from the regularity condition \eqref{eq:Tumarkin_condition}, the following example shows that it would be reasonable to have a representing \emph{distribution} instead of a representing \emph{measure}. Such kinds of representations may be found \eg in \cite{LangerWoracek2015}.

\begin{example}\label{ex:distributions}
Consider the function $q(z) := z^{-2}$ for $z \in \C\setminus\R$. This function fails to satisfy the regularity condition as
\begin{multline*}
\int_\R|q(t+\I\,r)-q(t-\I\,r)|\frac{\diff t}{1+t^2} = \int_\R\frac{4\,r|t|}{(1+t^2)(t^2+r^2)^2}\diff t \\
= \frac{4\,(1-r^2+2\,r^2\ln(r))}{r(1-r^2)^2} \xrightarrow{r \to 0^+} +\infty.
\end{multline*}
As such, it is not a rational quasi-Herglotz function and does not admit an integral representation of the form \eqref{eq:intRep_1var_quasi}. However, we note this particular function admits a representation of an analogous form as \eqref{eq:intRep_1var_quasi} where the measure $\nu$ is replaced by a distribution $V$ with compact support. Indeed, for any such distribution $V$ and $a,b \in \C$, the expression
$$q(z) = a + b\,z + \frac{1}{\pi}\big\langle V, t \mapsto \til{K}(z,t) \big\rangle$$
determines a well-defined function on $\C\setminus\R$. For the particular choice $a = -1$, $b = 0$ and $V = \pi\delta_0'$, \ie the derivative of the Dirac distribution, it holds that
$$-1 + \frac{1}{\pi}\big\langle \pi\delta_0', t \mapsto \til{K}(z,t) \big\rangle = \frac{1}{z^2}$$
for every $z \in \C\setminus\R$. Furthermore, for this function $q$, it still holds that
$$a = \frac{1}{2}\big(q(\I) + q(-\I)\big)$$
and that
$$\lim\limits_{z \ntto \infty}\frac{q(z)}{z} = \lim\limits_{z \nntto \infty}\frac{q(z)}{z} = b.$$
Note, however, that we do not discuss if such a representation is unique, {nor the most general class of distributions that could be used}.\hfill$\lozenge$
\end{example}

\section{Connections with other topics}\label{sec:comparisons}

\subsection{Weighted Hardy space $\hardy^1$}\label{subsec:Hardy}
For quasi-Herglotz functions that are identically zero in $\C^-$ the regularity condition \eqref{eq:Tumarkin_condition} reminds on the defining (but stronger) condition for a certain weighted Hardy space $\hardy^1(\C^+;w)$. In this section we are going to compare these spaces.

Let $f \in \hol(\C^+)$ and  let $w(x) := (1+x^2)^{-1}$ be a weight-function on $\R$. Consider the expression 
\begin{equation*}
    \label{eq:weighted_Hardy_norm}
    \|f\|_{1,w} := \sup\limits_{y > 0}\int_\R|f(x+\I\:y)|w(x)\diff x = \sup\limits_{y > 0}\int_\R|f(x+\I\:y)|\frac{\diff x}{1+x^2}.
\end{equation*}
The weighted Hardy space $\hardy^1(\C^+;w)$ is then defined as
$$\hardy^1(\C^+;w) := \{f \in \hol(\C^+)~|~\|f\|_{1,w} < \infty\}.$$

For any function $f \in \hardy^1(\C^+;w)$ we may now consider the function $q$ on $\C\setminus\R$, defined via
\begin{equation}\label{eq:hardy_temp}
q(z) := \left\{\begin{array}{rcl}
f(z) & ; & z \in \C^+, \\[0.2cm]
0 & ; & z \in \C^-.
\end{array}\right.
\end{equation}
This function surely satisfies the regularity condition \eqref{eq:Tumarkin_condition} as the requirement that $\|f\|_{1,w} < \infty$ is much stronger than condition \eqref{eq:Tumarkin_condition}. To show that the function $q$ also satisfies growth condition \eqref{eq:Vladimirov_growth}, it suffices to show this for the function $f$. To that end, consider the function $g\colon \C^+ \to \C$ given as
$$g(z) := \frac{f(z)}{(z+\I)^2}.$$
This function lies in the un-weighted Hardy space $\hardy^1(\C^+)$ as
$$\sup\limits_{y > 0}\int_\R|g(x+\I\,y)|\diff x = \sup\limits_{y > 0}\int_\R\frac{|f(x+\I\:y)|}{x^2+(1+y)^2}\diff x \leq \|f\|_{1,w} < \infty.$$
Due to a standard result for Hardy spaces, \cf \cite[Thm. 5.19]{RosenblumRovnyak1994}, the function $g$ has a boundary value $\til{g}$ on $\R$ almost everywhere. Furthermore, $\til{g} \in \mathrm{L}^1(\R)$ and it holds that
\begin{equation}
    \label{eq:Hardy_rep}
    g(z) = \frac{1}{2\pi\,\I}\int_\R\frac{\til{g}(t)}{t-z}\diff x
\end{equation}
for all $z = x+\I\,y \in \C^+$. Hence,
$$|g(x+\I\,y)| \leq \frac{1}{2\pi}\,\frac{\|\til{g}\|_{\mathrm{L}^1(\R)}}{y},$$
yielding, in terms of the function $f$, that
$$ |f(x+\I\,y)| \leq \frac{1}{2\pi}\,\|\til{g}\|_{\mathrm{L}^1(\R)}\,\frac{|x+\I\,y+\I|^2}{y} \leq M\,\frac{1+x^2+y^2}{y}$$
for some constant $M \geq 0$. Therefore, every function $f \in \hardy^1(\C^+;w)$ that is extended to a function $q$ {on} $\C\setminus\R$ via formula \eqref{eq:hardy_temp} gives rise to a quasi-Herglotz function. Note also that, for a function $g \in \hardy^1(\C^+)$, representation \eqref{eq:Hardy_rep} is automatically well-defined for all $z \in \C\setminus\R$ and it holds that $g(z) = 0$ for all $z \in \C^-$ \cite[Thm. 5.19]{RosenblumRovnyak1994}. Hence, in this sense $ \hardy^1(\C^+;w)$ is contained in the set of  quasi-Herglotz functions vanishing in the lower half plane. However, theses sets do not coincide. Namely, let $q$ be a quasi-Herglotz function that is identically zero in $\C^-$ and consider the function $f$ on $\C^+$ given as the restriction of the function $q$ to the upper half-plane, \ie
$$f(z) := q|_{\C^+}(z).$$
For such a function $f$, it does not necessarily hold that $\|f\|_{1,w} < \infty$. Indeed, consider the function 
$$q(z) := \left\{\begin{array}{rcl}
\sqrt{z} & ; & z \in \C^+, \\[0.2cm]
0 & ; & z \in \C^-,
\end{array}\right.$$
where the branch cut of the square-root is taken along the negative real axis. This function obviously satisfies the growth condition \eqref{eq:Vladimirov_growth} and we can easily check that is also satisfies the regularity condition \eqref{eq:Tumarkin_condition}, as it holds, for any $y \in [0,1]$, that
$$\int_\R|q(x+\I\:y)|\frac{\diff x}{1+x^2} = \int_\R\frac{(x^2+y^2)^{\frac{1}{4}}}{1+x^2}\diff x \leq \int_\R(1+x^2)^{-\frac{3}{4}}\diff x = \sqrt{\pi}\frac{\Gamma(\frac{1}{4})}{\Gamma(\frac{3}{4})},$$
where $\Gamma$ denotes Euler's Gamma function.

Let now $f$ be the restriction of $q$ to $\C^+$, \ie $f(z) = \sqrt{z}$ for $z \in \C^+$. For any $y > 0$, we estimate that
\begin{multline*}
\int_\R|f(x+\I\:y)|\frac{\diff x}{1+x^2} = \int_\R\frac{(x^2+y^2)^{\frac{1}{4}}}{1+x^2}\diff x = \sqrt{y}\int_\R\frac{(1+\frac{x^2}{y^2})^{\frac{1}{4}}}{1+x^2}\diff x \\
\geq \sqrt{y}\int_\R\frac{\diff x}{1+x^2} = \pi\sqrt{y} \xrightarrow{y \to \infty} \infty.
\end{multline*}
Hence, we have found that there exist quasi-Herglotz functions that are identically zero in $\C^-$ such that their restrictions to $\C^+$ do not belong to $\hardy^1(\C^+;w)$.

\subsection{Definitizable functions}\label{subsec:definitizable}

We also want to mention another class of functions that has non-empty intersection with quasi-Herglotz functions, namely the so-called definitizable functions, see \cite{Jonas2000,Jonas2006}.

\begin{define}\label{def:definitizable}
A function $g:\C\setminus\R\to\C$ is called \emph{definitizable} if it is piecewise meromorphic in $\C\setminus\R$, symmetric with respect to the real line, \ie $g(\overline z)=\overline{g(z)}$, and satisfies the following {three} conditions:
\begin{itemize}
\item[{\bf A:}]{the function $g$ has no more than a finite number of nonreal poles,}
\item[{\bf B:}]{the order of growth of $g$ near $\R$ is finite, \ie there exit constants $M$ and $m$ such that $|g(z)|\leq M\cdot\frac{(|z|+1)^m}{|\Im[z]|^m}$ for all $z$ in some neighbourhood of the closed real line $\overline{\R}$, and}
\item[{\bf C:}]{there is a finite (possibly empty) subset $E$ of $\R$ such that every connected component of $\R\setminus E$ is of definite (\ie either of positive or negative) type for $g$.}
\end{itemize}
Here, an open set $\Delta\subseteq\R$ is said to be \emph{of positive type for $g$} if for every sequence $(z_n)_{n\in\N}\subset\C^+$ which converges to a point in $\Delta$ it holds
\begin{equation}\label{eq:positive_type}
 \liminf_{n\to\infty}\Im[
 g(z_n)] \geq 0
\end{equation}
and, correspondingly, is said to be of \emph{negative type for $g$} if it is of positive type for {the function} $-g$. 
\end{define}

Alternatively, definitizable functions may also be characterized in the following way, \cf \cite[Def. 1.1]{Jonas2000}. A piecewise meromorphic  function $G$ in $\C\setminus\R$, symmetric with respect to $\mathbb R$ is definitizable if and only if there exists an $\mathbb R$-symmetric rational function $r$ such that the
product $rG$ is the sum of a {ordinary Herglotz} function $N$ and a rational
function $P$ with the poles of $P$ being points of holomorphy of $G$:
$$r(z)G(z) = N(z) + P(z)$$
for all points $z \in \C\setminus\R$ of holomorphy of $rG$. 

On one side, definitizable functions allow for higher order singularities than quasi-Herglotz functions (both in $\R$ and $\C\setminus\R$), but on the other side, the sign condition \textbf{C} is quite restrictive.

Without going into details, we would also like to mention that functions that are analytic in some subset of $\C^+$ can be represented with resolvents of self-adjoint operators (or, more generally, self-adjoint linear relations) in Kre\u{\i}n spaces (theses are certain vector spaces with an indefinite inner product), \cf \cite{DijksmaLangerdeSnoo1987}. For example, ordinary Herglotz functions correspond, in this sense, to self-adjoint linear relations in Hilbert spaces, \ie the inner product is positive definite and, hence, the spectrum (which in this case is confined to the real line) is of so-called positive type.

Using this fact shows directly that a real quasi-Herglotz function $q=h_1-h_2$ admits a corresponding representation in a Kre\u{\i}n space $\mathcal K=\mathcal H_1[+](-\mathcal H_2)$, and where the corresponding operator $A$ is in block diagonal form with respect to this fundamental decomposition of $\mathcal K$. 

For definitizable functions, however, this is not necessarily the case. These functions have, instead, representations with so-called definitizable operators. A bounded self-adjoint operator $A$ in a Kre\u{\i}n space $\mathcal K$ is called \emph{definitizable} if there exists a polynomial $p$ such that $[p(A)\mathbf x, \mathbf x]_{\mathcal K}\geq0$ for all $\mathbf x\in\mathcal K$, where $[\:\cdot\:,\:\cdot\:]_{\mathcal{K}}$ denotes the indefinite inner-product on $\mathcal{K}$. Hence, for definitzable functions, there are finitely many intervals of the real line, which are of positive type, and finitely many of negative type (\ie where $p(t)>0$ and $p(t)<0$, respectively.) Furthermore, at the zeros of $p$, the definitizable function may have higher order singularities.

For quasi-Herglotz functions, however, this is not true, \ie the sign type may change arbitrarily many times.  More precisely, let a real quasi-Herglotz function $q$ be given by its integral representation with a measure $\nu$, which has the Hahn decomposition $\nu=\nu_+-\nu_-$ into a difference of positive measures $\nu_\pm$. Then, in general, there does not exist a collection of \emph{finitely many} intervals such that the support of $\nu_+$ is contained in the closure of these intervals and the support of $\nu_-$ in the closure of the complement. 

The sets of quasi-Herglotz functions and definitizable functions do have a large intersection, \eg all ordinary Herglotz functions. However, none is contained in the other. Examples of definitizable functions that are not quasi-Herglotz functions are $g(z) = z^{-2}$ (see also Example \ref{ex:distributions}) and
$$g(z) = \left\{\begin{array}{rcl}
\I\,z & ; & z \in \C^+, \\
-\I\,z & ; & z \in \C^-,
\end{array}\right.$$
both of which fail to satisfy the regularity condition \eqref{eq:Tumarkin_condition}. Conversely, for the real quasi-Herglotz function
$$q(z) := \left\{\begin{array}{rcl}
\E^{\I\,z} & ; & z \in \C^+,  \\
\E^{-\I\,z} & ; & z \in \C^-, 
\end{array}\right.$$
given by the data $(\E^{-1},0,\nu)$, with $\diff\nu(t) := \sin(t)(1+t^2)^{-1}\diff t$, in the sense of Theorem \ref{thm:intRep_1var_quasi}, it holds that
$$\supp(\nu_+) = \bigcup_{k \in \Z}[2k\pi,(2k+1)\pi] \quad\text{and}\quad \supp(\nu_-) = \bigcup_{k \in \Z}[(2k-1)\pi,2k\pi]$$
and, hence, $q$ cannot be a definitizable function.

Finally, we would also like to mention that the well studied class of \emph{Generalized Nevanlinna functions}, \cite{KreinLanger1977,LangerWoracek2015}, which are special definitizable functions, has an intersection with the class of quasi-Herglotz functions, however, none of the two is contained in the other.

\subsection{Cauchy transform on the unit disk}\label{subsec:Cauchy}

The integral representation theorem \ref{thm:intRep_1var_quasi} can also be reinterpreted as the answer to the question which function on $\C\setminus\R$ can appear as the integral transform of a complex Borel measure $\nu$, where the transform is given by the kernel $\til{K}$. On the unit disk, a classical answer to an analogous question is known, but the integral kernel used there is not a direct biholomorphic transform of the kernel $\til{K}$. The classical setting on the unit circle is the following, \cf \cite{CimaMathesonRoss2006}.

Let $\sigma$ be a complex Borel measure on the unit circle $S^1$ and let $\hat{\C}$ denote the Riemann sphere, \ie $\hat{\C} = \C \cup \{\infty\}$ equipped with the standard topology of the sphere $S^2$. The \emph{Cauchy transform} $\mathfrak{C}$ of a measure $\sigma$ is a function on $\hat{\C}\setminus S^1$ defined as
$$(\mathfrak{C}\sigma)(\tau) := \int_{S^1}\frac{1}{1-\bar{\zeta}\,\tau}\diff\sigma(\zeta) = \int_{S^1}\frac{\zeta}{\zeta - \tau}\diff\sigma(\zeta).$$
Any measure $\sigma$ on $S^1$ can be transformed to a measure $\til{\sigma}$ on $[0,2\pi)$ via the change of variables $\zeta = \E^{\I\,s}$. In the $s$-variable, the Cauchy transform takes the form
\begin{equation}
    \label{eq:Cauchy_trans_theta}
    (\mathfrak{C}\sigma)(\tau) = \int_{[0,2\pi)}\frac{\I\,\E^{\I\,s}}{1-\E^{-\I\,s}\tau}\diff\til{\sigma}(s) = \frac{\I\,\til{\sigma}(\{0\})}{1-\tau} + \I\int_{(0,2\pi)}\frac{\E^{\I\,s}}{1-\E^{-\I\,s}\tau}\diff\til{\sigma}(s).
\end{equation}
Note that $(\mathfrak{C}\sigma)(\infty) = 0$ for any measure $\sigma$.

The Cayley transform $\psi$ can, in the present context, be viewed as an automorphism of $\hat{\C}$. It is given by the formula
$$\psi\colon \zeta \mapsto \I\,\frac{1+\zeta}{1-\zeta}$$
and maps
$$\psi\colon\left\{\begin{array}{LCL}
\D & \to & \C^+, \\
\C\setminus \bar{\D} & \to & \C^-\setminus\{-\I\}, \\
\{\infty\} & \to & \{-\I\}, \\
S^1\setminus\{1\} & \to & \R, \\
\{1\} & \to & \{\infty\}.
\end{array}\right.$$
As such, for $s \in (0,2\pi)$, we have $\psi(\E^{\I\,s}) \in \R$. The inverse Cayley transform $\varphi$ is given by
$$\varphi\colon \xi \mapsto \frac{\xi - \I}{\xi + \I}.$$
In particular, for $t \in \R$, it holds that $\varphi(t) \in S^1\setminus\{1\}$. If a change of variables between $s \in (0,2\pi)$ and $t \in \R$ is given by $\E^{\I\,s} = \varphi(t)$, then
$$\diff s = \frac{2}{1+t^2}\diff t.$$

A classical theorem of Tumarkin now characterizes which holomorphic functions appear as the Cauchy transform of a complex Borel measure on the unit circle, \cf \cite[Thm. 5.3.1]{CimaMathesonRoss2006} and \cite[Thm. 1]{Tumarkin1956}.

\begin{thm}[{Tumarkin}]
\label{thm:Tumarkhin_disk}
Let $F$ be a holomorphic function on $\hat{\C}\setminus S^1$ with $F(\infty) = 0$. Then, there exists a complex Borel measure $\sigma$ on $S^1$ such that $\mathfrak{C}\sigma = F$ if and only if
\begin{equation}
\label{eq:Tumarkhin_disk}
\sup\limits_{r \in (0,1)}\int_{S^1}|F(r\,\zeta) - F(r^{-1}\zeta)|\zeta^{-1}\diff\zeta < \infty.    
\end{equation}
\end{thm}

\begin{remark}
A generalization of this theorem to arbitrary domains appears in \cite{MarkushevichTumarkin1997} where the idea is to preserve the form of the integral representation, \ie have a representation of the form
$$\int_\gamma\frac{\diff\sigma(\zeta)}{\zeta - \tau},$$
while weakening the regularity requirement. Therefore, this generalization does not relate to quasi-Herglotz function in any stronger form than Theorem \ref{thm:Tumarkhin_disk} already does. Representations of the same form have also been considered in \eg \cite{Havin1958,Havin1963}, see also \cite{CimaMathesonRoss2006} for a general overview.
\end{remark}

Via the Cayley transform and its inverse, the information from Theorem \ref{thm:Tumarkhin_disk} can be translated to the case of quasi-Herglotz functions. If $q$ is a holomorphic function on $\C\setminus\R$, then $F = q \circ \psi$ is a holomorphic function on $\hat{\C}\setminus S^1$. The requirement that $F(\infty) = 0$ implies that we must have $q(-\I) = 0$. Note that this cannot be satisfied by any ordinary Herglotz function.

\begin{figure}[!ht]
\begin{tikzpicture}
\def\bigradius{2}
\def\shift{6.8}

\draw [help lines,->] (-1.25*\bigradius, 0) -- (1.25*\bigradius,0) node[above] {$x$};
\draw [help lines,->] (0, -1.25*\bigradius) -- (0, 1.25*\bigradius) node[right] {$\I\,y$};
\draw [dashed,line width=1pt] (0.9*\bigradius,0) arc (0:360:0.9*\bigradius);
\draw [line width=0.9pt,decoration={markings,
mark=at position 0.6 with {\arrow[line width=1pt]{>}},
mark=at position 0.88 with {\arrow[line width=1pt]{<}}
},postaction={decorate}] (0,0) -- (-1.1*\bigradius,0.66*\bigradius);
\fill (0,0)  circle[radius=1.5pt];
\node at (0.2,0.2) {$0$};
\fill (-0.8575*0.9*\bigradius,0.5145*0.9*\bigradius) circle[radius=1.5pt];
\fill [white] (-0.8575*0.9*\bigradius,0.5145*0.9*\bigradius) circle[radius=0.75pt];
\node at (-0.8575*0.9*\bigradius-0.08,0.5145*0.9*\bigradius+0.27) {$\zeta$};
\fill (-1.1*\bigradius,0.66*\bigradius) circle[radius=1.5pt];
\node at (-1.1*\bigradius-0.1,0.66*\bigradius+0.2) {$\infty$};

\draw [help lines,->] (-1.25*\bigradius+\shift, 0) -- (1.25*\bigradius+\shift,0) node[above] {$x$};
\draw [help lines,->] (0+\shift, -1.25*\bigradius) -- (0+\shift, 1.25*\bigradius) node[right] {$\I\,y$};
\draw [dashed,line width=1pt] (-1.1*\bigradius+\shift,0) -- (1.1*\bigradius+\shift,0);
\draw [line width=0.9pt,decoration={markings,
mark=at position 0.45 with {\arrow[line width=1pt]{>}},
mark=at position 0.59 with {\arrow[line width=1pt]{<}}
},postaction={decorate}] ([shift=(58:1.18)]{-0.6222+\shift},0) arc (58:302:1.18);
\fill (0+\shift,\bigradius/2) circle[radius=1.5pt];
\node at (0.25+\shift,\bigradius/2) {$+\I$};
\fill (0+\shift,-\bigradius/2) circle[radius=1.5pt];
\node at (0.25+\shift,-\bigradius/2) {$-\I$};
\fill (-0.9*\bigradius+\shift,0) circle[radius=1.5pt];
\fill [white] (-0.9*\bigradius+\shift,0)  circle[radius=0.75pt];
\node at (-0.9*\bigradius+\shift-0.15,0.2) {$t$};
\end{tikzpicture}
    
\caption{Approaching a point $\zeta \in S^1$ in condition \eqref{eq:Tumarkhin_disk} and the corresponding situation on the real line.}
\label{fig:regularity}
\end{figure}
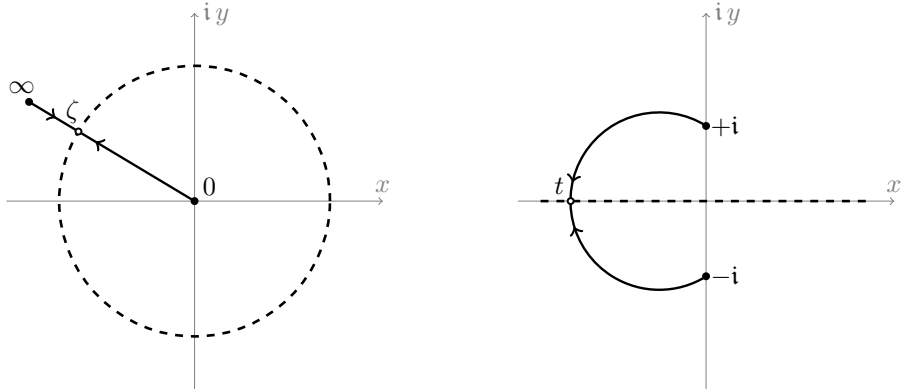

Furthermore, the regularity condition \eqref{eq:Tumarkhin_disk} is related to the regularity condition \eqref{eq:Tumarkin_condition}. Condition \eqref{eq:Tumarkhin_disk} can first be rewritten as
$$\sup\limits_{r \in (0,1)}\int_{[0,2\pi)}|q(\psi(r\,\E^{\I\,s})) - q(\psi(r^{-1}\,\E^{\I\,s})|\diff s < \infty,$$
where a factor of $\I$ from the change of variables was thrown away as it does not influence the finiteness of the above supremum. Since this integral is weighted against the Lebesgue measure on $[0,2\pi)$, we may skip integration over the point at zero. Using the bijection between $(0,2\pi)$ and $\R$ mentioned earlier, we may further rewrite condition \eqref{eq:Tumarkhin_disk} as
$$\sup\limits_{r \in (0,1)}\int_{\R}|q(\psi(r\,\varphi(t))) - q(\psi(r^{-1}\,\varphi(t)))|\frac{\diff t}{1+t^2} < \infty,$$
where a factor of $2$ was thrown away as before. In condition \eqref{eq:Tumarkhin_disk}, for a fixed $\zeta \in S^1$, the functions $r \mapsto r\,\zeta$ and $r \mapsto r^{-1}\,\zeta$ parametrize, for $r \in (0,1)$, the two segments of the radial line between $0$ and $\infty$ passing through $\zeta$. After the transformation, for a fixed $t \in \R$, the function $r \mapsto \psi(r\,\varphi(t))$ parametrizes, for $r \in (0,1)$, the circular arc between $+\I$ and $t$ that approaches the real line at the angle $\frac{\pi}{2}$ while the function $r \mapsto \psi(r^{-1}\,\varphi(t))$ parametrizes the mirrored arc between $-\I$ and $t$. This is visualized in Figure \ref{fig:regularity}.

\begin{figure}[!ht]
\begin{tikzpicture}
\def\bigradius{2}
\def\shift{6.8}

\draw [help lines,->] (-1.25*\bigradius, 0) -- (1.25*\bigradius,0) node[above] {$x$};
\draw [help lines,->] (0, -1.25*\bigradius) -- (0, 1.25*\bigradius) node[right] {$\I\,y$};
\draw [dashed,line width=1pt] (-1.1*\bigradius,0) -- (1.1*\bigradius,0);
\draw [help lines] (-1.25*\bigradius,1.4) -- (1.25*\bigradius,1.4);
\draw [help lines] (-1.25*\bigradius,-1.4) -- (1.25*\bigradius,-1.4);
\draw [line width=0.9pt,decoration={markings,
mark=at position 0.42 with {\arrow[line width=1pt]{>}},
mark=at position 0.65 with {\arrow[line width=1pt]{<}}
},postaction={decorate}] (1.5,-1.4) -- (1.5,1.4);
\fill (1.5,-1.4) circle[radius=1.5pt];
\node at (1.5+0.4,-1.4-0.2) {$t-\I$};
\fill (1.5,0) circle[radius=1.5pt];
\fill [white] (1.5,0) circle[radius=0.75pt];
\node at (1.5+0.2,0+0.2) {$t$};
\fill (1.5,1.4) circle[radius=1.5pt];
\node at (1.5+0.4,1.4+0.2) {$t+\I$};

\draw [help lines,->] (-1.25*\bigradius+\shift, 0) -- (1.25*\bigradius+\shift,0) node[above] {$x$};
\draw [help lines,->] (0+\shift, -1.25*\bigradius) -- (0+\shift, 1.25*\bigradius) node[right] {$\I\,y$};
\draw [dashed,line width=1pt] (0.9*\bigradius+\shift,0) arc (0:360:0.9*\bigradius);
\draw [help lines] (0.9*\bigradius+\shift,-1.25*\bigradius) -- (0.9*\bigradius+\shift,1.25*\bigradius);
\draw [help lines] (0.9*\bigradius+\shift,0) arc (0:360:0.45*\bigradius);
\draw [line width=0.9pt,decoration={markings,
mark=at position 0.52 with {\arrow[line width=1pt]{>}},
mark=at position 0.79 with {\arrow[line width=1pt]{<}}
},postaction={decorate}] ([shift=(90:1.03923)]{1.8+\shift},1.03923) arc (90:190:1.03923);
\fill (1.8+\shift,2.07846) circle[radius=1.5pt];
\node at (1.8+\shift+0.3,2.07846) {$\zeta_-$};
\fill (0.9+\shift,1.559) circle[radius=1.5pt];
\fill [white] (0.9+\shift,1.559) circle[radius=0.75pt];
\node at (0.9+\shift-0.15,1.559+0.3) {$\zeta$};
\fill (0.771429+\shift,0.890769) circle[radius=1.5pt];
\node at (0.771429+\shift,0.890769-0.25) {$\zeta_+$};
\end{tikzpicture}
    
\caption{Approaching a point $t \in \R$ in condition \eqref{eq:Tumarkin_condition} and the corresponding situation on the unit circle.}
\label{fig:regularity_v2}
\end{figure}
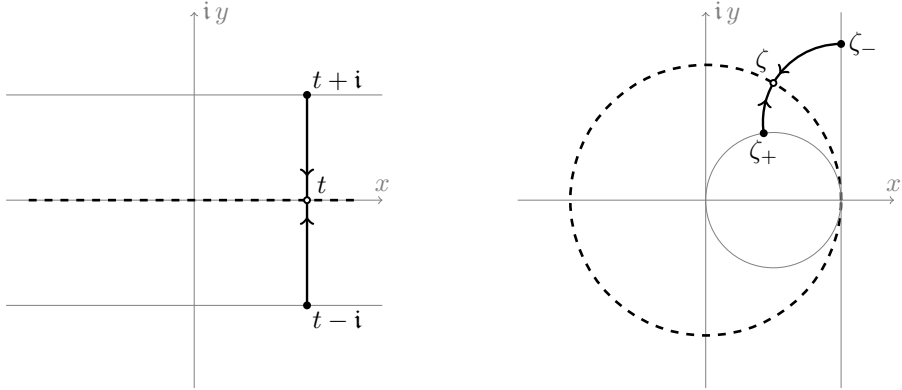

Conversely, in the regularity condition \eqref{eq:Tumarkin_condition}, we approach a given number $t \in \R$ via the straight lines between $t \pm \I$ and $t$. Under the inverse Cayley transform $\varphi$, the point $t$ maps to a point $\zeta \in S^1\setminus\{1\}$, while
$$\varphi(t+\I) = \frac{1+\zeta}{3-\zeta} =: \zeta_+  \quad\text{and}\quad \varphi(t-\I) = -\left(\frac{1+\zeta}{3-\zeta}\right)^{-1} =: \zeta_- .$$
The circle containing the points $\zeta,\zeta_+$ and $\zeta_-$ always intersects the unit circle at an angle of $\frac{\pi}{2}$. Moreover, the point $\zeta_+ \in \D$ always lies on the circle $\{|z-\frac{1}{2}|=\frac{1}{2}\}$ while the point $\zeta_- \in \C\setminus\bar{\D}$ always lies on the circle $\{\Re[z] = 1\}\cup\{\infty\}$. This is visualized in Figure \ref{fig:regularity_v2}.

The following theorem now establishes the precise relation between Cauchy transforms on the unit circle and quasi-Herglotz functions.

\begin{thm}\label{thm:Cauchy}
Let $q\colon \C\setminus\R \to \C$ be a holomorphic function. Then, $q$ is a quasi-Herglotz function if and only if there exists a number $c \in \C$ and a complex Borel measure $\sigma$ on $S^1$ such that it holds, for every $z \in \C\setminus\R$, that
\begin{equation}
    \label{eq:Cauchy_q}
    q(z) = c + (\mathfrak{C}\sigma)(\varphi(z)),
\end{equation}
where $\varphi$ denotes the inverse Cayley transform and $\mathfrak{C}$ denotes the Cauchy transform. Furthermore, the number $c$ is given in terms of the function $q$ (or its data $(a,b,\nu)$) as
\begin{equation}
    \label{eq:Cauchy_a}
    c = q(-\I) = a - \I\left(b + \frac{1}{\pi}\int_\R\diff\nu(t)\right).
\end{equation}
\end{thm}

\proof
Assume first that we have a function $q$ on $\C\setminus\R$ defined via equality \eqref{eq:Cauchy_q} for some number $c$ and some measure $\sigma$ as in the theorem. Using the second form in equality \eqref{eq:Cauchy_trans_theta} to write the Cauchy transform $\mathfrak{C}\sigma$ of the measure $\sigma$, we infer via equality \eqref{eq:Cauchy_q} that the function $q$ admits, for every $z \in \C\setminus\R$, an integral representation of the form
$$q(z) = c + \frac{\I\,\til{\sigma}(\{0\})}{1-\varphi(z)} + \I\int_{(0,2\pi)}\frac{\E^{\I\,s}}{1-\E^{-\I\,s}\varphi(z)}\diff\til{\sigma}(s),$$
where the measure $\sigma$ on $S^1$ has been reparametrized to a measure $\til{\sigma}$ on $[0,2\pi)$ as before. Changing the variable $s \in (0,2\pi)$ to $t \in \R$ by setting  $\E^{\I\,s} = \varphi(t)$ yields
$$q(z) = c + \frac{\I\,\til{\sigma}(\{0\})}{1-\varphi(z)} + \I\int_{\R}\frac{\varphi(t)}{1-\frac{\varphi(z)}{\varphi(t)}}\cdot\frac{2}{1+t^2}\diff\til{\mu}(t),$$
where $\til{\mu}$ is the transform of the measure $\til{\sigma}$ under this change of variables.
Simplifying the above expression now yields
\begin{multline}
\label{eq:q_cauchy_rep}
q(z) = c + \frac{z+\I}{2}\,\til{\sigma}(\{0\}) + \int_\R\frac{(t-\I)(z+\I)}{(t-z)(t+\I)^2}\diff\til{\mu}(t) \\ = c +  \frac{z+\I}{2}\,\beta + \int_\R\frac{(t-\I)^2(z+\I)}{(t-z)(t+\I)}\diff\til{\nu}(t),
\end{multline}
where $\beta := \til{\sigma}(\{0\})$ and
$$\diff\til{\nu}(t) := \frac{1}{1+t^2}\diff\til{\mu}(t).$$

We claim now that representation \eqref{eq:q_cauchy_rep} is actually of the form \eqref{eq:intRep_1var_quasi}, implying, by Theorem \ref{thm:intRep_1var_quasi}, that the function $q$ is, in fact, a quasi-Herglotz function. To do that, we will show how the data $(a,b,\nu)$ may be defined in terms of the numbers $c,\beta$ and the measure $\til{\nu}$, yielding the desired result by the uniqueness statement of Theorem \ref{thm:intRep_1var_quasi}.

As such, we calculate that
\begin{equation}
    \label{eq:Cauchy_a_temp}
    \frac{1}{2}\big(q(\I)+q(-\I)\big) = c + \I\left(\frac{\beta}{2} + \int_\R\frac{t-\I}{t+\I}\diff\til{\nu}(t)\right)
\end{equation}
and
$$\lim\limits_{z \ntto \infty}\frac{q(z)}{z} = \lim\limits_{z \nntto \infty}\frac{q(z)}{z} = \frac{\beta}{2},$$
while it holds that, for any function $g$ as in Proposition \ref{prop:Stieltjes_inversion}, that
\begin{multline*}
    \lim\limits_{y \to 0^+}\int_\R g(x)\tfrac{1}{2\I}(q(x+\I\,y)-q(x-\I\,y))\diff x \\
    = \lim\limits_{y \to 0^+}\int_\R g(x)\left(\frac{\beta\,y}{2(1+t^2)} + \int_\R\pois(x+\I\,y,t)(t-\I)^2\diff\til{\nu}(t)\right)\diff x \\
    = \pi\int_\R g(t)(t-\I)^2\diff\til{\nu}(t) = \pi\int_\R g(t)(1+t^2)\frac{t-\I}{t+\I}\diff\til{\nu}(t).
\end{multline*}
Here, we used the Fubini-Tonelli theorem to change the order of integration with respect to the $t$- and $x$-variables. Note also that the $\pi$-factor comes from the Poisson kernel. 
If we now define
$$a := c + \I\left(\frac{\beta}{2} + \int_\R\frac{t-\I}{t+\I}\diff\til{\nu}(t)\right), \quad b := \frac{\beta}{2}$$
and
$$\diff\nu(t) := \pi\frac{t-\I}{t+\I}\diff\til{\nu}(t),$$
then representation \eqref{eq:intRep_1var_quasi} becomes representation \eqref{eq:q_cauchy_rep}, as desired. Furthermore, returning with this information to equality \eqref{eq:Cauchy_a_temp}, we infer that the representing parameters $a$, $b$ and $\nu$ for the function $q$ and the number $c$ do indeed satisfy equality \eqref{eq:Cauchy_a}.

Conversely, assume that we start with a quasi-Herglotz function $q$. Then, one may define a number $c \in \C$ using equality \eqref{eq:Cauchy_a} and the integral representation \eqref{eq:intRep_1var_quasi} for the function $q$ may be rewritten as
\begin{multline*}
    q(z) = c + \I\left(b + \frac{1}{\pi}\int_\R\diff\nu(t)\right) + b\,z + \frac{1}{\pi}\int_\R \til{K}(z,t)\diff\nu(t) \\
    = c + b(z+\I) + \frac{1}{\pi}\int_\R\left(\til{K}(z,t)+\I\right)\diff\nu(t) \\
    = c + b(z+\I) + \frac{1}{\pi}\int_\R\frac{(t-\I)(z+\I)}{(t-z)(t+\I)^2}\cdot(t+\I)^2\diff\nu(t).
\end{multline*}
Define now $\beta := 2b$ and
$$\diff\til{\mu}(t) := \frac{1}{\pi}\,(t+\I)^2\diff\nu(t).$$
Using this, we may now construct a measure $\til{\sigma}$ on $[0,2\pi)$ by setting $\til{\sigma}(\{0\}) := \beta = 2b$ and choosing $\til{\sigma}|_{(0,2\pi)}$ to be determined by the change of variables $\E^{\I\,s} = \varphi(t)$, \ie
$$\diff\til{\sigma}|_{(0,2\pi)}(s) = \frac{2}{1+t^2}\diff\til{\mu}(t).$$
The measure $\til{\sigma}$ can then be reparametrized to a measure $\sigma$ on $S^1$. By reversing the calculations made before, equality \eqref{eq:Cauchy_q} holds for this particular measure $\sigma$ and the number $c$ as defined before, finishing the proof. 
\endproof

Theorem \ref{thm:Cauchy} shows that Cauchy transforms of complex Borel measures on $S^1$ form a strict subclass of quasi-Herglotz functions when mapped over to $\C\setminus\R$ via the Cayley transform. In particular, this subclass does not include any non-trivial ordinary Herglotz function, but does include all quasi-Herglotz functions that are identically zero in the lower half-plane. However, the full class of quasi-Herglotz functions may be recovered with the addition of a complex constant.

\subsection{Sum-rules}\label{subsec:sumrules}

As mentioned in the introduction, one motivation to study quasi-Herglotz functions comes from applications, see \cite{IvanenkoETAL2019b}, and is also related to so-called sum-rules.  Let us briefly recall what is meant by these. 

Ordinary Herglotz functions admit asymptotic expansions  with real coefficients around any point $t_0 \in \R$ or at the point at infinity as long as the expansion is restricted to some (upper or lower) Stoltz domain. More precisely, let $h$ be an ordinary Herglotz function. Then, we say that $h$ admits, at $z=t_0$, an \emph{asymptotic expansion of order $M\geq-1$} if there exist real numbers $a_{-1},a_0,a_1,\ldots,a_M$ (depending on $t_0$) such that $h$ can be written as
\begin{equation}\label{eq:asymp_zero}
\displaystyle{h(z) = \frac{a_{-1}}{z-t_0} + a_{0} + a_1(z-t_0) + \ldots + a_M(z-t_0)^M + o\big((z-t_0)^M\big)}
\end{equation}
as $z \ntto t_0$ or $z \nntto t_0$. Similarly, we say that $h$ admits, at $z=\infty$, an \emph{asymptotic expansion of order $K\geq-1$} if there exist real numbers $b_1,b_0,b_{-1} ,\ldots,b_{-K}$ such that $h$  can be written as
\begin{equation}\label{eq:asymp_infty}
\displaystyle{h(z) = b_{1} z + b_{0} +\frac{b_{-1}}z+\ldots + \frac{b_{-K}}{z^{K}}+ o\Big(\frac{1}{z^{K}}\Big)}
\end{equation}
as $z \ntto \infty$ or $z \nntto \infty$.

Indeed, at $z = 0$, an expansion of order $M = -1$ always exists for any Herglotz function $h$ due to Proposition \ref{prop:quasi_point_mass_limits}, see also \eg \cite{BernlandEtal2011,KacKrein1974}. Similarly, at $z = \infty$, an expansion of order $K = -1$ due to formula \eqref{eq:b_constant}. Furthermore, the number $b$ from representation \eqref{eq:intRep_1var} and formula \eqref{eq:b_constant} equals the number $b_1$ appearing in expansion \eqref{eq:asymp_infty}.

For an ordinary Herglotz function, there is a close relation between the above asymptotic expansions and  certain weighted integrals. More precisely, an ordinary Herglotz function $h$ admits an asymptotic expansion  at $\infty$ of order $2N_\infty\geq0$ if and only if the following limit exists 
\begin{equation}\label{eq:sum_rule_limit}
\lim\limits_{\varepsilon\to 0^+}\lim\limits_{y \to 0^+}\int_{\varepsilon<|x|<\frac{1}{\varepsilon}}x^{2N_\infty}\Im[h(x+\I\,y)]\diff x.
\end{equation}
In this case,  for $k=0,1,2,\ldots, 2N_\infty$ it holds 
\begin{equation*}\label{eq:sum_rule_identities}
\lim\limits_{\varepsilon\to 0^+}\lim\limits_{y \to 0^+}\int_{\varepsilon<|x|<\frac{1}{\varepsilon}}x^{k}\Im[h(x+\I\,y)]\diff x \\ = 
a_{-k-1} - b_{-k-1},
\end{equation*}
where $a_{-k-1}$ and $b_{-k-1}$ are the coefficients from the expansions \eqref{eq:asymp_zero} and \eqref{eq:asymp_infty}, respectively. Note that only in the case $k=0$ two terms appear, otherwise, for positive $k$, the $a$-coefficient vanishes. A similar relation holds for weight functions with negative exponents, which are then related to the expansion at $z=0$. These identities are known as \emph{sum-rules} and are \eg used for \emph{a priori} estimates for physical limitations in passive systems \cite{BernlandEtal2011}. 

Already for real quasi-Herglotz function, an analogous statement does not necessarily hold. Recently, \cf \cite[Thm. 3.3]{IvanenkoETAL2019b}, it was noted that an analogous result holds for a function $q = h_1 - h_2$, with $h_1,h_2$ ordinary Herglotz functions, as long as one assumes \emph{a priori} that at least one of the functions $h_1,h_2$ admits an asymptotic expansion of some order. Contrary to ordinary Herglotz functions, a real quasi-Herglotz functions may admit \eg a series expansion of the form \eqref{eq:asymp_infty} of order $K = 0$ (\ie a higher order than the guaranteed one), but the left-hand side of the corresponding formula \eqref{eq:sum_rule_identities} does not exist as a finite number. An example of such a function is $q(z) := \tan(z) - \I$, which admits, at $z = \infty$, an asymptotic expansion of the form
$$q(z) = o(1),$$
\ie of order $2N_\infty$ with $N_\infty = 0$. The corresponding identity \eqref{eq:sum_rule_identities}, for $k = 0$, would say that
$$\lim\limits_{\varepsilon\to 0^+}\lim\limits_{y \to 0^+}\int_{\varepsilon<|x|<\frac{1}{\varepsilon}}\Im[q(x+\I\,y)]\diff x = a_{-1} - b_{-1} = 0.$$
However,
\begin{multline*}\int_{\varepsilon<|x|<\frac{1}{\varepsilon}}\Im[q(x+\I\,y)]\diff x 
= \int_{\varepsilon<|x|<\frac{1}{\varepsilon}}\left(\frac{\sh(2y)}{\cos(2x)+\ch(2y)}-1\right)\diff x \\
= 2\left(\arctan(\tan(1/\varepsilon)\th(y))-\arctan(\tan(\varepsilon)\th(y)) + \varepsilon - \frac{1}{\varepsilon}\right) \\
\xrightarrow{y \to 0^+} 2\left(\varepsilon - \frac{1}{\varepsilon}\right) \xrightarrow{\varepsilon \to 0^+} -\infty.
\end{multline*}

\subsection*{Acknowledgement}
The authors thank Harald Woracek for his help concerning Hardy-space functions.

\bibliographystyle{amsplain}
\bibliography{total}

\end{document}